\documentclass[reqno,11pt]{amsart} 
\usepackage{amsmath}
\usepackage{amssymb}
\usepackage{amsthm}
\usepackage{verbatim}
\usepackage{xcolor}
\usepackage{graphics}

\newcommand\NoBlackBoxes{\global\overfullrule0pt}
\NoBlackBoxes
\theoremstyle{plain} 

\def\4{\kern1pt}

\def\6{\vphantom0}

\def\8{\kern-10pt}
\def\7#1{_{(#1)}}

\makeatletter
\let\serieslogo@\relax
\let\@setcopyright\relax

\def\speciallabelmark#1{\def\@currentlabel{#1}}
\makeatother

\begin{document}

\def\ffrac#1#2{\raise.5pt\hbox{\small$\4\displaystyle\frac{\,#1\,}{\,#2\,}\4$}}
\def\ovln#1{\,{\overline{\!#1}}}
\def\ve{\varepsilon}
\def\kar{\beta_r}

\title{REGULARIZED DISTRIBUTIONS AND ENTROPIC \\ 
STABILITY OF CRAMER'S CHARACTERIZATION \\
OF THE NORMAL LAW
}

\author{S. G. Bobkov$^{1,4}$}
\thanks{1) School of Mathematics, University of Minnesota, USA;
Email: bobkov@math.umn.edu}
\address
{Sergey G. Bobkov \newline
School of Mathematics, University of Minnesota  \newline 
127 Vincent Hall, 206 Church St. S.E., Minneapolis, MN 55455 USA
\smallskip}
\email {bobkov@math.umn.edu} 

\author{G. P. Chistyakov$^{2,4}$}
\thanks{2) Faculty of Mathematics, University of Bielefeld, Germany;
Email: chistyak@math.uni-bielefeld.de}
\address
{Gennadiy P. Chistyakov\newline
Fakult\"at f\"ur Mathematik, Universit\"at Bielefeld\newline
Postfach 100131, 33501 Bielefeld, Germany}
\email {chistyak@math.uni-bielefeld.de}

\author{F. G\"otze$^{3,4}$}
\thanks{3) Faculty of Mathematics, University of Bielefeld, Germany;
Email: goetze@math.uni-bielefeld.de}
\thanks{4) Research partially supported by NSF grant DMS-1106530
and SFB 701}
\address
{Friedrich G\"otze\newline
Fakult\"at f\"ur Mathematik, Universit\"at Bielefeld\newline
Postfach 100131, 33501 Bielefeld, Germany}
\email {goetze@mathematik.uni-bielefeld.de}



\begin{abstract}
For regularized distributions we establish stability of the 
characterization of the normal law in Cramer's theorem with respect 
to the total variation norm and the entropic distance.
As part of the argument, Sapogov-type theorems are refined for
random variables with finite second moment.
\end{abstract}

\maketitle
\markboth{S. G. Bobkov, G. P. Chistyakov and F. G\"otze}{Stability
in Cramer's theorem}

\def\theequation{\thesection.\arabic{equation}}
\def\E{{\bf E}}
\def\R{{\bf R}}
\def\C{{\bf C}}
\def\P{{\bf P}}
\def\H{{\rm H}}
\def\Im{{\rm Im}}
\def\Tr{{\rm Tr}}

\def\k{{\kappa}}
\def\M{{\cal M}}
\def\Var{{\rm Var}}
\def\Ent{{\rm Ent}}
\def\O{{\rm Osc}_\mu}

\def\ep{\varepsilon}
\def\phi{\varphi}
\def\F{{\cal F}}
\def\L{{\cal L}}

\def\be{\begin{equation}}
\def\en{\end{equation}}
\def\bee{\begin{eqnarray*}}
\def\ene{\end{eqnarray*}}


\section{{\bf Introduction}}
\setcounter{equation}{0}

\vskip2mm
Let $X$ and $Y$ be independent random variables. A theorem of Cramer [Cr] 
indicates that, if the sum $X+Y$ has a normal distribution, then both 
$X$ and $Y$ are normal. P. L\'evy established stability of this 
characterization property with respect to the L\'evy distance, which 
is formulated as follows. Given $\ep > 0$ and distribution functions $F$, $G$,
$$
L(F*G, \Phi) < \ep \ \ \Rightarrow \ \
L(F,\Phi_{a_1,\sigma_1}) < \delta_\ep, \
L(G,\Phi_{a_2,\sigma_2}) < \delta_\ep,
$$
for some $a_1,a_2 \in \R$ and $\sigma_1, \sigma_2 > 0$, where $\delta_\ep$ only 
depends on $\ep$, and in a such way that $\delta_\ep \rightarrow 0$ as 
$\ep \rightarrow 0$. Here $\Phi_{a,\sigma}$ stands for the distribution 
function of the normal law $N(a,\sigma^2)$ with mean $a$ and standard deviation 
$\sigma$, i.e., with density
$$
\varphi_{a,\sigma}(x) = \frac{1}{\sigma\sqrt{2\pi}}\,e^{-(x-a)^2/2\sigma^2},
\quad x\in \R,
$$  
and we omit indices in the standard case $a=0$, $\sigma=1$. As usual,
$F*G$ denotes the convolution of the corresponding distributions.

The problem of quantitative versions of this stability property of the 
normal law has been intensively studied in many papers, starting 
with results by Sapogov [S1-3] and ending with results by Chistyakov 
and Golinskii [C-G], who found the correct asymptotic of the best possible
error function $\ep \rightarrow \delta_\ep$ for the L\'evy distance. 
See also [Z1], [M], [L-O], [C], [Se] [Sh1-2].

As for stronger metrics, not much is known up to now. According to 
McKean ([MC], cf. also [C-S] for some related aspects 
of the problem), it was Kac who raised the question about the stability 
in Cramer's theorem with respect to the entropic distance to normality.
Let us recall that, if a random variable $X$ with finite second moment 
has a density $p(x)$, its entropy
$$
h(X) = -\int_{-\infty}^\infty p(x) \log p(x)\,dx
$$
is well-defined and is bounded from above by the entropy of the normal 
random variable $Z$, having the same variance $\sigma^2 = \Var(Z) = \Var(X)$. 
The entropic distance to the normal is given by the formula
$$
D(X) = h(Z) - h(X) = \int_{-\infty}^\infty 
p(x) \log \frac{p(x)}{\varphi_{a,\sigma}(x)}\,dx,
$$
where in the last formula it is assumed that $a = \E Z = \E X$.
It represents the Kullback-Leibler distance from the distribution 
$F$ of $X$ to the family of all normal laws on the line.

In general, $0 \leq D(X) \leq \infty$, and an infinite value is 
possible. This quantity does not depend on the 
variance of $X$ and is stronger than the total variation distance 
$\|F - \Phi_{a,\sigma}\|_{{\rm TV}}$, as may be seen from the Pinsker 
(Pinsker-Csisz\'ar-Kullback) inequality
$$
D(X)\, \geq\, \frac{1}{2}\, \|F - \Phi_{a,\sigma}\|_{{\rm TV}}^2.
$$

Thus, Kac's question is whether one can bound the
entropic distance $D(X+Y)$ from below in terms of $D(X)$ and $D(Y)$
for independent random variables, i.e., to have an inequality
$$
D(X+Y) \geq \alpha(D(X),D(Y))
$$
with some non-negative function $\alpha$, such that $\alpha(t,s) > 0$ for $t,s>0$. 
If so, Cramer's theorem would be an immediate consequence of this. Note that the 
reverse inequality does exist, and in case $\Var(X+Y)=1$ we have
\bee
D(X+Y) \leq \Var(X) D(X) + \Var(Y) D(Y),
\ene
which is due to the general entropy power inequality, cf. [D-C-T].

It turned out that Kac's question has a negative solution. More 
precisely, for any $\ep > 0$, one can construct independent random 
variables $X$ and $Y$ with absolutely continuous symmetric distributions 
$F$, $G$, and with $\Var(X)=\Var(Y)=1$, such that

\vskip2mm
$a)$ \ $D(X+Y) < \ep\,;$

\vskip2mm
$b)$ \ $\|F - \Phi_{a,\sigma}\|_{{\rm TV}} > c$ \ and \ 
$\|G - \Phi_{a,\sigma}\|_{{\rm TV}} > c$, \
for all\, $a \in \R$ and $\sigma > 0$, 

\vskip2mm
\noindent
where $c > 0$ is an absolute constant, see [B-C-G1]. In particular,
$D(X)$ and $D(Y)$ are bounded away from zero. Moreover, refined analytic 
tools show that the random variables may be chosen to be identically
distributed, i.e., $a)-b)$ hold with $F = G$, see [B-C-G2].

Nevertheless, Kac's problem remains to be of interest for subclasses
of probability measures obtained by convolution with a ``smooth''
distribution. The main purpose of this note is to give an affirmative solution 
to the problem in the (rather typical) situation, when 
independent Gaussian noise is added to the given random variables. 
That is, for a small parameter $\sigma > 0$, we consider the 
regularized random variables
$$
X_\sigma = X + \sigma Z, \qquad Y_\sigma = Y + \sigma Z,
$$
where $Z$ denotes a standard normal random variable, independent of
$X,Y$. As a main result, we prove:

\vskip5mm
{\bf Theorem 1.1.} {\it Let $X,Y$ be independent random variables 
with $\Var(X+Y) = 1$. Given $0 < \sigma \leq 1$,
the regularized random variables $X_\sigma$ and $Y_\sigma$ satisfy
$$
D(X_\sigma + Y_\sigma) \geq 
\exp\Big\{- \frac{c\, \log^7(2+1/D)}{D^2}\Big\},
$$
where $c>0$ is an absolute constant, and
$$
D = \sigma^2\, 
\big(\Var(X_\sigma)\, D(X_\sigma) + \Var(Y_\sigma)\, D(Y_\sigma)\big).
$$
}

Thus, if $D(X_\sigma + Y_\sigma)$ is small, the entropic distances 
$D(X_\sigma)$ and $D(Y_\sigma)$ have to be small, as well. 
In particular, Cramer's theorem is a consequence of this statement.
However, it is not clear whether the above lower bound is optimal with 
respect to the couple $(D(X_\sigma),D(Y_\sigma))$, and perhaps the
logarithmic term in the exponent may be removed. As we will see, 
a certain improvement of the bound can be achieved, when $X$ and $Y$ have equal 
variances.

Beyond the realm of results around P. L\'evy's theorem, recently there has been 
renewed the interest in other related stability problems in different areas of 
Analysis and Geometry. One can mention, for example, the problems of sharpness 
of the Brunn-Minkowski and Sobolev-type inequalities (cf. [F-M-P1-2], [Seg], [B-G-R-S]).

%
%

We start with the description and refinement of Sapogov-type theorems 
about the normal approximation in Kolmogorov distance (Sections 2-3) 
and then turn to analogous results for the L\'evy distance (Section 4). 
A version of Theorem 1.1 for the total variation distance is given in Section 5. 
Sections 6-7 deal with the problem of bounding the tail function
$\E X^2\, 1_{\{|X|\geq T\}}$ in terms of the entropic distances
$D(X)$ and $D(X+Y)$, which is an essential part of Kac's problem. 
A first application, namely, to a variant of Chistyakov-Golinskii's theorem, 
is discussed in Section 8. In Section 9, we develop several estimates connecting 
the entropic distance $D(X)$ and the uniform deviation of the density $p$ from 
the corresponding normal density. 
In Section 10 an improved variant of Theorem 1.1 is derived
in the case, where $X$ and $Y$ have equal variances. 
The general case is treated in Section 11.
Finally, some relations between different distances in the space
of probability distributions on the line are postponed to appendix.

%
%
%
%
%
%
%
%
%
%
%
%
%
%
%


\vskip5mm
\section{{\bf Sapogov-type theorems for Kolmogorov distance}}
\setcounter{equation}{0}

\vskip2mm
Throughout the paper we consider the following classical metrics 
in the space of probability distributions on the real line:

1) The Kolmogorov or $L^\infty$-distance
$\|F - G\| = \sup_x {|F(x) - G(x)|}$;

2) The L\'evy distance
$$
L(F,G) = \min\big\{h \geq 0: \
G(x-h)-h \leq F(x) \leq G(x+h)+h, \ \ \forall x \in \R\big\};
$$

3) The Kantorovich or $L^1$-distance
$$
W_1(F,G) = \int_{-\infty}^\infty |F(x) - G(x)|\,dx;
$$

4) The total variation distance
$$
\|F - G\|_{{\rm TV}} = 
\sup \sum |(F(x_k) - G(x_k)) - (F(y_k) - G(y_k))|,
$$
where the sup is taken over all finite collections of points
$y_1 < x_1 < \dots < y_n < x_n$.

In these relations, $F$ and $G$ are arbitrary distribution functions.
Note that the quantity $W_1(F,G)$ is finite, as long as both $F$ and 
$G$ have a finite first absolute moment.


In the sequel, $\Phi_{a,v}$ or $N(a,v^2)$ denote the normal 
distribution (function) with parameters $(a,v^2)$, $a \in \R$ $v>0$. 
If $a = 0$, we write $\Phi_v$, and write $\Phi$ 
in the standard case $a=0$, $v = 1$.

Now, let $X$ and $Y$ be independent random variables with distribution
functions $F$ and $G$. Then the convolution $F*G$ represents 
the distribution of the sum $X+Y$.
If both random variables have mean zero and unit variances,
Sapogov's main stability result reads as follows: 

\vskip4mm
{\bf Theorem 2.1.} {\it Suppose that $\E X = \E Y = 0$ and 
$\Var(X) = \Var(Y) = 1$. If
$$
\|F * G - \Phi * \Phi\| \leq \ep < 1,
$$
then with some absolute constant $C$
$$
\|F - \Phi\| \leq \frac{C}{\sqrt{\log\frac{1}{\ep}}} \qquad and \qquad
\|G - \Phi\| \leq \frac{C}{\sqrt{\log\frac{1}{\ep}}}\, .
$$
}

In the general case (that is, when there are no finite moments), the 
conclusion is somewhat weaker. Namely, with $\ep \in (0,1)$, we 
associate
$$
a_1 = \int_{-N}^N x\,dF(x), \quad
\sigma_1^2 = \int_{-N}^N x^2\,dF(x) - a_1^2 \ \ \ (\sigma_1 \geq 0),
$$
and similarly $(a_2,\sigma_2^2)$ for the distribution function $G$, 
where $N = N(\ep) = 1 + \sqrt{2\log(1/\ep)}$.

In the sequel, we also use the function
$$
m(\sigma,\ep) = 
\min\Big\{\frac{1}{\sqrt{\sigma}},\log \log \frac{e^e}{\ep}\Big\}, \qquad
\sigma>0, \ \ 0 < \ep \leq 1.
$$

\vskip4mm
{\bf Theorem 2.2.} {\it Assume $\|F * G - \Phi\| \leq \ep < 1$.
If $F$ has median zero, and $\sigma_1,\sigma_2 > 0$, then with some 
absolute constant $C$
$$
\|F - \Phi_{a_1,\sigma_1}\|\, \leq\, 
\frac{C}{\sigma_1 \sqrt{\log\frac{1}{\ep}}}\,m(\sigma_1,\ep),
$$
and similarly for $G$.
}

\vskip4mm
Originally, Sapogov derived a weaker bound in [Sa1-2] with worse 
behaviour with respect to both $\sigma_1$ and $\ep$. In [Sa3] he gave 
an improvement,
$$
\|F - \Phi_{a_1,\sigma_1}\|\, \leq\, 
\frac{C}{\sigma_1^3 \, \sqrt{\log\frac{1}{\ep}}}
$$
with a correct asymptotic of the right-hand side with respect to $\ep$,
cf. also [L-O]. The correctness of the asymptotic with respect to $\ep$
was studied in [M], cf. also [C]. In 1976 Senatov [Se1], using 
the ridge property of characteristic functions, improved the factor
$\sigma_1^3$ to $\sigma_1^{3/2}$, i.e., 
\be
\|F - \Phi_{a_1,\sigma_1}\|\, \leq\, 
\frac{C}{\sigma_1^{3/2} \, \sqrt{\log\frac{1}{\ep}}}.
\en
He also emphasized that the presence of $\sigma_1$ in the bound is 
essential. A further improvement of the power of $\sigma_1$ is due to 
Shiganov [Sh1-2]. Moreover, at the expense of an additional $\ep$-dependent
factor, one can replace $\sigma_1^{3/2}$ with $\sigma_1$.
As shown in [C-G], see Remark on p.\,2861,
\be
\|F - \Phi_{a_1,\sigma_1}\|\, \leq\, 
\frac{C\log \log \frac{e^e}{\ep}}{\sigma_1 \sqrt{\log\frac{1}{\ep}}}.
\en
Therefore, Theorem 2.2 is just the combination of the two results,
(2.1) and (2.2).

Let us emphasize that all proofs of these theorems use
the methods of the Complex Analysis. Moreover, up to now there 
is no "Real Analysis" proof of the Cram\'er theorem and of its extensions
in the form of Sapogov-type results. This, however, does not concern
the case of identically distributed summands, cf. [B-C-G2].

We will discuss the bounds in the L\'evy distance in the next sections.

The assumption about the median in Theorem 2.2 may be weakened to the
condition that the medians of $X$ and $Y$, $m(X)$ and $m(Y)$, are 
bounded in absolute value by a constant. For example, if 
$\E X = \E Y = 0$ and $\Var(X+Y) = 1$, and if, for definiteness, 
$\Var(X) \leq 1/2$, then, by Chebyshev's inequality, $|m(X)| \leq 1$, 
while $|m(Y)|$ will be bounded by an absolute constant, 
when $\ep$ is small enough, due to the main hypothesis
$\|F*G - \Phi\| \leq \ep$.

Moreover, if the variances of $X$ and $Y$ are bounded away from zero, 
the statement of Theorem 2.2 holds with $a_1 = 0$, and the factor 
$\sigma_1$ can be replaced with the standard deviation of $X$. 
In the next section, we recall some standard arguments in order 
to justify this conclusion and give a more general version of 
Theorem 2.2 involving variances:

\vskip4mm
{\bf Theorem 2.3.} {\it Let $\E X = \E Y = 0$, $\Var(X+Y) = 1$. If\, 
$\|F * G - \Phi\| \leq \ep < 1$, then with some absolute constant $C$
$$
\|F - \Phi_{v_1}\|\, \leq\, 
\frac{Cm(v_1,\ep)}{v_1 \sqrt{\log\frac{1}{\ep}}} \quad and \quad
\|G - \Phi_{v_2}\|\, \leq\, 
\frac{Cm(v_2,\ep)}{v_2 \sqrt{\log\frac{1}{\ep}}},
$$
where $v_1^2 = \Var(X)$, $v_2^2 = \Var(Y)$ $(v_1,v_2 > 0)$.
}

\vskip4mm
Under the stated assumptions, Theorem 2.3 is stronger than Theorem 2.2, 
since $v_1 \geq \sigma_1$. Another advantage of this formulation is 
that $v_1$ does not depend on $\ep$, while $\sigma_1$ does.


\vskip5mm
\section{{\bf Proof of Theorem 2.3}}
\setcounter{equation}{0}

\vskip2mm
Let $X$ and $Y$ be independent random variables with distribution 
functions $F$ and $G$, respectively, with $\E X = \E Y = 0$ and 
$\Var(X+Y) = 1$. We assume that
$$
\|F * G - \Phi\| \leq \ep < 1,
$$
and keep the same notations as in Section 2.
Recall that $N = N(\ep) = 1 + \sqrt{2 \log(1/\ep)}$.

The proof of Theorem 2.3 is entirely based on Theorem 2.2. We will need: 

\vskip4mm
{\bf Lemma 3.1.} {\it With some absolute constant $C$ we have
$$
0 \leq 1 - (\sigma_1^2 + \sigma_2^2) \leq CN^2\sqrt{\ep}.
$$
}

A similar assertion, $|\sigma_1^2 + \sigma_2^2 - 1| \leq CN^2 \ep$, is known 
under the assumption that $F$ has a median at zero (without moment assumptions). 
For the proof of Lemma 3.1, we use arguments from [Sa1] and [Se1], cf. Lemma~1. 
It will be convenient to divide the proof into several steps.

\vskip4mm
{\bf Lemma 3.2.} {\it Let\ 
$\ep \leq \ep_0 = \frac{1}{4} - \Phi(-1) = 0.0913...$
Then $|m(X)| \leq 2$ and $|m(Y)| \leq 2$.
}

\vskip5mm
Indeed, let $\Var(X) \leq 1/2$. Then $|m(X)| \leq 1$, 
by Chebyshev's inequality. Hence,
$$
\frac{1}{4} \leq \P\{X \leq 1, \, Y \leq m(Y)\} \leq 
\P\{X + Y \leq m(Y)+1\} \leq \Phi(m(Y)+1) + \ep,
$$
which for $\ep \leq \frac{1}{4}$ implies that 
$m(Y)+1 \geq \Phi^{-1}(\frac{1}{4} - \ep)$. In particular,
$m(Y) \geq -2$, if $\ep \leq \ep_0$. Similarly, $m(Y) \leq 2$.
\qed

To continue, introduce truncated random variables at level $N$. Put
$X^* = X$ in case $|X| \leq N$, $X^* = 0$ in case $|X| > N$, and
similarly $Y^*$ for $Y$. Note that 
$$
\E X^* = a_1, \ \ \Var(X^*) = \sigma_1^2, \quad {\rm and} \quad
\E Y^* = a_2, \ \ \Var(Y^*) = \sigma_2^2.
$$
By the construction, $\sigma_1 \leq v_1$ and $\sigma_2 \leq v_2$. 
In particular,
$
\sigma_1^2 + \sigma_2^2 \leq v_1^2 + v_2^2 = 1.
$
Let $F^*,G^*$ denote the distribution functions of $X^*,Y^*$,
respectively. 

\vskip4mm
{\bf Lemma 3.3.} {\it  With some absolute constant $C$ we have
$$
\|F^* - F\| \leq C\sqrt{\ep}, \qquad \|G^* - G\| \leq C\sqrt{\ep}, \qquad
\|F^* * G^* - \Phi\| \leq C\sqrt{\ep}.
$$
}

{\bf Proof.}
One may assume that $N = N(\ep)$ is a point of continuity of both $F$ and $G$. 
Since the Kolmogorov distance is bounded by 1, one may also assume that $\ep$ 
is sufficiently small, e.g., $\ep < \min\{\ep_0,\ep_1\}$, where
$\ep_1 = \exp\{-1/(3 - 2\sqrt{2}\,)\}$. In this case 
$(N-2)^2 > (N-1)^2/2$, so
$$
\Phi(-(N-2)) = 1 - \Phi(N-2) \leq \frac{1}{2}\,e^{-(N-2)^2/2} \leq
\frac{1}{2}\,e^{-(N-1)^2/4} = \frac{\sqrt{\ep}}{2}.
$$
By Lemma 3.2 and the basic assumption on the convolution $F * G$,
\bee
\frac{1}{2}\,\P\{Y \leq -N\}
  &\leq&
\P\{X \leq 2, \, Y \leq -N\} \\
 & \leq & 
\P\{X + Y \leq -(N-2)\}  =  (F * G)(-(N-2))
 \, \leq \,  \Phi(-(N-2)) + \ep.
\ene
So, $G(-N) \leq 2\Phi(-(N-2)) + 2\ep \leq 3 \sqrt{\ep}$. Analogously,
$1-G(N) \leq 3 \sqrt{\ep}$. Thus,
\bee
\int_{\{|x| \geq N\}} dG(x) \leq 6 \sqrt{\ep} \qquad
\text{as well as}\qquad
\int_{\{|x| \geq N\}} dF(x) \leq 6 \sqrt{\ep}.
\ene

In particular, for $x < -N$, we have $|F^*(x) - F(x)| = F(x) \leq 6 \sqrt{\ep}$,
and similarly for $x > N$. If $-N < x < 0$, then $F^*(x) = F(x)-F(-N)$, and if 
$0 < x < N$, we have $F^*(x) = F(x)+(1-F(N))$. In both cases, 
$|F^*(x) - F(x)| \leq 6 \sqrt{\ep}$. Therefore,
$$
\|F^* - F\| \leq 6 \sqrt{\ep}.
$$
Similarly, $\|G^* - G\| \leq 6 \sqrt{\ep}$. From this, by the triangle inequality,
\bee
\|F^* * G^* - F * G\| & \leq & 
\|F^* * G^* - F^* * G\| + \|F^* * G - F * G\| \\
 & \leq &
\|F^* - F\| + \|G^* - G\| \, \leq \, 12 \sqrt{\ep}.
\ene
Finally, 
$$
\|F^* * G^* - \Phi\| \leq \|F^* * G^* - F * G\| + \|F * G - \Phi\| \leq 
12\, \sqrt{\ep} + \ep \leq 13\, \sqrt{\ep}.
$$
\qed

{\bf Proof of Lemma 3.1.} 
Since $|X^* + Y^*| \leq 2N$ and
$a_1 + a_2 = \E\, (X^* + Y^*) = \int x\, dF^* * G^*(x)$, we have,
integrating by parts,
\bee
a_1 + a_2 & = & \int_{-2N}^{2N} x\, d ((F^* * G^*)(x) - \Phi(x)) \\
 & = &
x\, ((F^* * G^*)(x) - \Phi(x))\bigg|_{x=-2N}^{x=2N} -
\int_{-2N}^{2N} ((F^* * G^*)(x) - \Phi(x))\,dx.
\ene
Hence, 
$
|a_1 + a_2| \leq 8N\,\|F^* * G^* - \Phi\|,
$
which, by Lemma 3.3, is bounded by $CN \sqrt{\ep}$.
Similarly,
\bee
\E\, (X^* + Y^*)^2 - 1 
 & = & 
\int_{-2N}^{2N} x^2\, d ((F^* * G^*)(x) - \Phi(x)) 
- \int_{\{|x|>2N\}} x^2\, d\Phi(x)
\\
 & = &
x^2\, ((F^* * G^*)(x) - \Phi(x))\bigg|_{x=-2N}^{x=2N} \\
 & &
 -\ 2 \int_{-2N}^{2N} x\, ((F^* * G^*)(x) - \Phi(x))\,dx
 - \int_{\{|x|>2N\}} x^2\, d\Phi(x).
\ene
Hence, 
$$
\left|\E\, (X^* + Y^*)^2 - 1\right|\, \leq\, 
24\,N^2\,\|F^* * G^* - \Phi\| +
2 \int_{2N}^\infty x^2\, d\Phi(x).
$$
The last integral asymptotically behaves like
$2N \varphi(2N) < N e^{-2(N-1)^2} = N \ep^4$. Therefore,
$\left|\E\, (X^* + Y^*)^2 - 1\right|$ is bounded by $CN^2 \sqrt{\ep}$.
Finally, writing
$\sigma_1^2 + \sigma_2^2 = \E\, (X^* + Y^*)^2 - (a_1 + a_2)^2$, 
we get that
$$
\left|\sigma_1^2 + \sigma_2^2 - 1\right| \leq 
\left|\E\, (X^* + Y^*)^2 - 1\right| + (a_1 + a_2)^2 \leq CN^2 \sqrt{\ep}
$$
with some absolute constant $C$. Lemma 3.1 follows.
\qed

\vskip4mm
{\bf Proof of Theorem 2.3.}
First note that, given $a > 0$, $\sigma > 0$, and $x \in \R$, the function
$$
\psi(x) = \Phi_{0,\sigma}(x) - \Phi_{a,\sigma}(x) =
\Phi\big(\frac{x}{\sigma}\big) - \Phi\big(\frac{x-a}{\sigma}\big)
$$
is vanishing at infinity, has a unique extreme point $x_0 = \frac{a}{2}$, and
$
\psi(x_0) = \int_{-a/2\sigma}^{a/2\sigma} \varphi(y)\,dy \leq
\frac{a}{\sigma\sqrt{2\pi}}.
$
Hence, including the case $a \leq 0$, as well, we get
$$
\|\Phi_{a,\sigma} - \Phi_{0,\sigma}\|\, \leq\, 
\frac{|a|}{\sigma \sqrt{2\pi}}.
$$

We apply this estimate for $a = a_1$ and $\sigma = \sigma_1$.
Since $\E X = 0$ and $\Var(X+Y) = 1$, by Cauchy's and Chebyshev's
inequalities,
$$
|a_1| = \left|\E\, X\, 1_{\{|X| \geq N\}}\right| \leq
\P\{|X| \geq N\}^{1/2} \leq \frac{1}{N} <
\frac{1}{\sqrt{\log\frac{e}{\ep}}}.
$$
Hence,
$$
\|\Phi_{a_1,\sigma_1} - \Phi_{0,\sigma_1}\|\, \leq\, 
\frac{|a_1|}{\sigma_1 \sqrt{2\pi}} \leq 
\frac{C}{\sigma_1 \sqrt{\log\frac{1}{\ep}}}.
$$
A similar inequality also holds for the parameters $(a_2,\sigma_2)$.

Now, define the non-negative numbers $u_1 = v_1 - \sigma_1$, 
$u_2 = v_2 - \sigma_2$. By Lemma 3.1,
\bee
CN^2 \sqrt{\ep} \ \geq \ 1 - (\sigma_1^2 + \sigma_2^2)
 & = &
1 - \left((v_1-u_1)^2 + (v_2-u_2)^2\right) \\
 & = & u_1\, (2v_1 - u_1) + u_2\, (2v_2 - u_2)
 \ \geq \ u_1 v_1 + u_2 v_2.
\ene
Hence,
$$
u_1 \leq \frac{CN^2 \sqrt{\ep}}{v_1} \quad {\rm and} \quad
u_2 \leq \frac{CN^2 \sqrt{\ep}}{v_2}.
$$
These relations can be used to estimate the Kolmogorov distance
$
\Delta = \|\Phi_{0,v_1} - \Phi_{0,\sigma_1}\|.
$

Given two parameters $\alpha > \beta > 0$, consider the function of 
the form $\psi(x) = \Phi(\alpha x) - \Phi(\beta x)$. In case $x>0$, 
by the mean value theorem, for some $x_0 \in (\beta x,\alpha x)$,
$$
\psi(x) = (\alpha - \beta)\, x \varphi(x_0) <
(\alpha - \beta)\, x \varphi(\beta x).
$$
Here, the right-hand side is maximized for $x = \frac{1}{\beta}$, which gives
$
\psi(x) < \frac{1}{\sqrt{2\pi e}}\, \frac{\alpha - \beta}{\beta}.
$
A similar bound also holds for $x<0$. Using this bound with 
$\alpha = 1/\sigma_1$ ($\sigma_1 > 0$), $\beta = 1/v_1$, we obtain
$$
\Delta \, \leq \,
\frac{1}{\sqrt{2\pi e}}\, v_1 \bigg(\frac{1}{\sigma_1} - \frac{1}{v_1}\bigg)
\, = \, \frac{1}{\sqrt{2\pi e}}\ \frac{u_1}{\sigma_1}
 \, \leq \, \frac{CN^2 \sqrt{\ep}}{\sigma_1 v_1}
 \, \leq \, \frac{CN^2 \sqrt{\ep}}{\sigma_1^2}.
$$

Thus, applying Theorem 2.2, we get with some 
universal constant $C>1$ that
\begin{eqnarray}
\|F - \Phi_{0,v_1}\| & \leq  &
\|F - \Phi_{a_1,\sigma_1}\| + 
\|\Phi_{a_1,\sigma_1} - \Phi_{0,\sigma_1}\| +
\|\Phi_{0,\sigma_1} - \Phi_{0,v_1}\| \nonumber \\
 & \leq  &
\frac{C}{\sigma_1 \sqrt{\log\frac{1}{\ep}}}\,m(\sigma_1,\ep) + 
\frac{C}{\sigma_1 \sqrt{\log\frac{1}{\ep}}} + 
\frac{CN^2 \sqrt{\ep}}{\sigma_1^2} \nonumber \\
 & \leq  &
\frac{2C}{\sigma_1 \sqrt{\log\frac{1}{\ep}}}\,m(\sigma_1,\ep) + 
\frac{CN^2 \sqrt{\ep}}{\sigma_1^2}. 
\end{eqnarray}
The obtained estimate remains valid when $\sigma_1 = 0$, as well.
On the other hand,
$
\sigma_1 = v_1 - u_1 \geq v_1 - \frac{CN^2 \sqrt{\ep}}{v_1} \geq 
\frac{1}{2}\, v_1
$
where the last inequality is fulfilled for the range
$v_1 \geq v(\ep) = \sqrt{C}\, N\, (4\ep)^{1/4}$. Hence, from (3.1)
and using $m(\sigma_1,\ep) \leq 2 m(v_1,\ep)$, for this range
$$
\|F - \Phi_{0,v_1}\| \, \leq \, 
\frac{8Cm(v_1,\ep)}{v_1 \sqrt{\log\frac{1}{\ep}}} + 
\frac{4CN^2 \sqrt{\ep}}{v_1^2}.
$$
Here, since $m(v_1,\ep) \geq 1$, the first term on the right-hand 
side majorizes the second one, if
$$
v_1 \geq \tilde v(\ep) = N^2\sqrt{\ep \log\frac{1}{\ep}}.
$$
Therefore, when $v_1 \geq w(\ep) = \max\{v(\ep),\tilde v(\ep)\}$,
with some absolute constant $C'$ we have
$$
\|F - \Phi_{0,v_1}\| \, \leq \, 
\frac{C'm(v_1,\ep)}{v_1 \sqrt{\log\frac{1}{\ep}}}.
$$
Thus, we arrive at the desired inequality for the range $v_1 \geq w(\ep)$. 
But the function $w$ behaves almost 
polynomially near zero and admits, for example, a bound of the form
$w(\ep) \leq \sqrt{C''}\, \ep^{1/6}$, $0 < \ep < \ep_0$,
with some universal $\ep_0 \in (0,1)$, $C''>1$.
So, when $v_1 \leq w(\ep)$, $0 < \ep < \ep_0$, we have
$$
\frac{1}{v_1 \sqrt{\log\frac{1}{\ep}}} \geq 
\frac{1}{w(\ep) \sqrt{\log\frac{1}{\ep}}} \geq
\frac{1}{\ep^{1/6} \sqrt{C''\log\frac{1}{\ep}}}.
$$
Here, the last expression is greated than 1, as long as $\ep$ is 
sufficiently small, say, for all $0 < \ep < \ep_1$, where $\ep_1$ is
determined by $(C'',\ep_0)$. Hence, for all such $\ep$, we have a better bound
$$
\|F - \Phi_{0,v_1}\| \, \leq \, \frac{C}{v_1 \sqrt{\log\frac{1}{\ep}}}.
$$

It remains to increase the constant $C'$ in order to involve the 
remaining values of $\ep$. A similar conclusion is true for the 
distribution $G$. Theorem 2.3 is thus proved completely.
\qed


\vskip5mm
\section{{\bf Stability in Cramer's theorem for the L\'evy distance}}
\setcounter{equation}{0}

\vskip2mm
Let $X$ and $Y$ be independent random variables with distribution functions $F$ and $G$. 
It turns out that in the bound of Theorem 2.2, the parameter
$\sigma_1$ can be completely removed, if we consider the stability 
problem for the L\'evy distance. More precisely, the following theorem 
was established in [C-G].

\vskip4mm
{\bf Theorem 4.1.} {\it Assume that\, $\|F * G - \Phi\| \leq \ep < 1$.
If $F$ has median zero, then with some absolute constant $C$
$$
L(F,\Phi_{a_1,\sigma_1})\, \leq\, C\,
\frac{(\log \log\frac{4}{\ep})^2}{\sqrt{\log\frac{1}{\ep}}}.
$$
}

Recall that
$$
a_1 = \int_{-N}^N x\,dF(x), \quad
\sigma_1^2 = \int_{-N}^N x^2\,dF(x) - a_1^2 \ \ \ (\sigma_1 \geq 0),
$$
and similarly $(a_2,\sigma_2^2)$ for the distribution function $G$, 
where $N = 1+\sqrt{2\log(1/\ep)}$.

As we have already discussed, the assumption about the median may be
relaxed to the condition that the median is bounded (by a universal 
constant).

The first quantitative stability result for the L\'evy distance, namely,
$$
L(F,\Phi_{a_1,\sigma_1})\, \leq\, C\, \log^{-1/8}(1/\ep),
$$
was obtained in 1968 by Zolotarev [Z1], who applied his famous 
Berry-Esseen-type bound. The power $1/8$ was later improved to $1/4$ 
by Senatov [Se1] and even more by Shiganov [Sh1-2]. The stated 
asymptotic in Theorem 4.1 is unimprovable, which was also shown in [C-G].

Note that in the assumption of Theorem 4.1, the Kolmogorov distance 
can be replaced with the L\'evy distance $L(F,\Phi)$ 
in view of the general relations
$$
L(F,\Phi) \leq \|F * G - \Phi\| \leq (1+M)\, L(F,\Phi), \qquad 
M = \|\Phi\|_{{\rm Lip}} = \frac{1}{\sqrt{2\pi}}.
$$
However, in the conclusion such replacement cannot be done at the 
expense of a universal constant, since we only have
$$
\|F - \Phi_{a_1,\sigma_1}\| \leq (1+M) \, L(F,\Phi_{a_1,\sigma_1}), \qquad 
M = \|\Phi_{a_1,\sigma_1}\|_{{\rm Lip}} = \frac{1}{\sigma_1\sqrt{2\pi}}.
$$

Now, our aim is to replace in Theorem 4.1 the parameters $(a_1,\sigma_1)$, 
which depend on $\ep$, with $(0,v_1)$ like in Theorem 2.3. 
That is, we have the following:

\vskip2mm
{\bf Question.} Assume that $\E X = \E Y = 0$, $\Var(X+Y) = 1$, and
$L(F * G,\Phi) \leq \ep < 1$. Is it true that
$$
L(F,\Phi_{v_1})\, \leq\, C\
\frac{(\log \log\frac{4}{\ep})^2}{\sqrt{\log\frac{1}{\ep}}}
$$
with some absolute constant $C$, where $v_1^2 = \Var(X)$?

\vskip2mm
In a sense, it is the question on the closeness of $\sigma_1$ to $v_1$
in the situation, where $\sigma_1$ is small. Indeed, using the triangle 
inequality, one can write
$$
L(F,\Phi_{v_1}) \leq
L(F,\Phi_{a_1,\sigma_1}) + L(\Phi_{a_1,\sigma_1},\Phi_{0,\sigma_1}) + 
L(\Phi_{\sigma_1},\Phi_{v_1}).
$$
Here, the first term may be estimated according to Theorem 4.1. For 
the second one, we have a trivial uniform bound (over all $\sigma_1$),
$$
L(\Phi_{a_1,\sigma_1},\Phi_{0,\sigma_1}) \leq |a_1|,
$$
which follows from the definition of the L\'evy metric.
In turn, the parameter $a_1$ admits the bound, which was already 
used in the proof of Theorem 2.3, namely, $|a_1| < \frac{1}{\sqrt{\log\frac{e}{\ep}}}$.
This bound behaves better than the one in Theorem 4.1, so we obtain:

\vskip4mm
{\bf Lemma 4.2.} {\it If\, $\E X = \E Y = 0$, $\Var(X+Y) = 1$, and
$L(F * G,\Phi) \leq \ep < 1$, then
$$
L(F,\Phi_{v_1}) \ \leq \
C\ \frac{(\log \log\frac{4}{\ep})^2}{\sqrt{\log\frac{1}{\ep}}} + 
L(\Phi_{\sigma_1},\Phi_{v_1}). \\
$$
}

Thus, we are reduced to estimating the distance $L(\Phi_{\sigma_1},\Phi_{v_1})$, 
which in fact should be done in terms of $v_1^2 - \sigma_1^2$.

\vskip4mm
{\bf Lemma 4.3.} {\it Given $v \geq \sigma \geq 0$, such that 
$v^2 - \sigma^2 \leq 1$, we have
$$
L(\Phi_{\sigma},\Phi_{v})^2 \leq (v^2 - \sigma^2)\, 
\log\frac{2}{v^2 - \sigma^2}.
$$
}


{\bf Proof.} It will be clear that the asymptotic in terms of 
$\alpha = \sqrt{v^2 - \sigma^2}$ is correct.

Since the normal distributions with mean zero are symmetric
about the origin, the L\'evy distance $L(\Phi_{\sigma},\Phi_{v})$ 
represents an optimal value $h \geq 0$, such that the inequality
\be
\Phi_{\sigma}(x) \leq \Phi_{v}(x+h) + h
\en
holds true for all $x$. (The other inequality, 
$\Phi_{v}(x) \leq \Phi_{\sigma}(x+h) + h$, is equivalent to (4.1)).
Moreover, for $x \leq 0$, we have 
$\Phi_{\sigma}(x) \leq \Phi_{v}(x)$, so only $x>0$ should be
taken into consideration. 

We may assume $v > \sigma$, i.e., $\alpha > 0$. Changing the variable
$x = \sigma y$, $y > 0$, (4.1) becomes
\be
\Phi(y) \leq \Phi\bigg(\frac{\sigma y+h}{v}\bigg) + h.
\en
Here $h$ needs to serve for all $\sigma>0$, while $\alpha$ is fixed. 
So, we need to minimize the function
$
\psi(\sigma) = \frac{\sigma y+h}{v} = 
\frac{\sigma y+h}{\sqrt{\sigma^2 + \alpha^2}}.
$
By the direct differentiation, we find that
$$
\psi'(\sigma) = \frac{y\alpha^2 - \sigma h}{(\sigma^2 + \alpha^2)^{3/2}} = 0 \
\Longleftrightarrow \ \sigma = \sigma_0 \equiv \frac{y\alpha^2}{h}, \qquad
\psi(\sigma_0) = \sqrt{y^2 + (h/\alpha)^2} \geq y.
$$
Since $\psi'(0) > 0$, we may conclude that $\psi$ is increasing for 
$\sigma  \leq \sigma_0$ and is decreasing for $\sigma  \geq \sigma_0$. 
Hence, the inequality (4.2) will only be strengthened, if we replace it with
\bee
\Phi(y) \ \leq \
\inf_{\sigma \geq 0} \Phi(\psi(\sigma)) + h \, = \,
\min\{\Phi(\psi(0)),\Phi(\psi(\infty))\} + h \, = \,
\min\Big\{\Phi\Big(\frac{h}{\alpha}\Big),\Phi(y)\Big\} + h.
\ene
That is,
$\Phi(y) \leq \Phi(\frac{h}{\alpha}) + h$, and since 
$y>0$ is arbitrary, it is equivalent to
$1 \leq \Phi(\frac{h}{\alpha}) + h$. In other words,
$$
L(\Phi_{\sigma},\Phi_{v}) \leq L(\Phi_0,\Phi_\alpha),
$$
where $\Phi_0$ denotes the unit mass at the origin. 

Thus, we are reduced to the case $\sigma = 0$. But then the L\'evy 
distance $h_0 = L(\Phi_0,\Phi_\alpha)$ represents the (unique) solution 
to the equation
$
1 = \Phi\big(\frac{h}{\alpha}\big) + h.
$
To estimate it, we may use the bound 
$1 - \Phi(\frac{h}{\alpha}) \leq \frac{1}{2}\, e^{-h^2/2\alpha^2}$
which gives $2h_0 \leq e^{-h_0^2/2\alpha^2}$. After the change 
$h_0 = \alpha \sqrt{2 \log(c_0/\alpha)}$, and using $\alpha \leq 1$, 
we obtain
$$
1 \geq 2c_0  \sqrt{2 \log\frac{c_0}{\alpha}} \geq 
2c_0  \sqrt{2 \log c_0},
$$
so, $4c_0^2\, \log c_0^2 \leq 1$. It follows that $c_0^2 < 2$ and
$h_0 \leq \alpha \sqrt{\log(2/\alpha^2)}$, as was claimed.
\qed

\vskip4mm
{\bf Remark.} Attempts to derive bounds on the L\'evy distance
$L(\Phi_{\sigma},\Phi_{v})$ by virtue of standard general relations, 
such as Zolotarev's Berry-Esseen-type estimate [Z2], lead to worse 
dependences of $\alpha^2 = v^2 - \sigma^2$.
For example, using a general relation $L(F,G)^2 \leq W_1(F,G)$, 
cf. Proposition A.1.1, together with the Kantorovich-Rubinshtein 
theorem, we get that
$$
L(\Phi_{\sigma},\Phi_{v})^2 \leq W_1(\Phi_{\sigma},\Phi_{v})
\leq \E\, |\sigma Z - v Z| \leq  v - \sigma =
\frac{v^2 - \sigma^2}{v + \sigma},
$$
where $Z \sim N(0,1)$ and where we did not loose much when bounding
$W_1$. This estimate has a disadvantage in comparison 
with Lemma 4.3 because of a possible small denominator.

\vskip4mm
In view of Lemmas 4.2-4.3, in order to proceed, one needs to bound 
$v_1^2 - \sigma_1^2$ in terms of $\ep$. However, this does not seem to be 
possible in general without stronger hypotheses. Note that
$$
v_1^2 - \sigma_1^2 = \int_{\{|x| > N\}} x^2\,dF(x) + a_1^2.
$$
Hence, we need to deal with the quadratic tail function
$$
\delta_X(T) = \int_{\{|x| > T\}} x^2\,dF(x), \qquad T \geq 0,
$$
whose behavior at infinity will play an important role in the sequel. 

Now, combining Lemmas 4.2 and 4.3, we obtain
$$
L(F,\Phi_{v_1}) \leq C\
\frac{(\log \log\frac{4}{\ep})^2}{\sqrt{\log\frac{1}{\ep}}} +
R\left(\delta_X(N) + a_1^2\right),
$$
where $R(t) = \sqrt{t \log(2/t)}$. This function is non-negative and 
concave in the interval $0 \leq t \leq 2$, with $R(0)=0$.
Hence, it is subadditive in the sense that
$R(\xi + \eta) \leq R(\xi) + R(\eta)$, for all $\xi, \eta \geq 0$,
$\xi + \eta \leq 2$. Hence,
\bee
R\left(\delta_X(N) + a_1^2\right) & \leq & R(\delta_X(N)) + R(a_1^2) \\
 & = & \bigg(\delta_X(N) \log\frac{2}{\delta_X(N)}\bigg)^{1/2} +
\sqrt{a_1^2 \log\frac{2}{a_1^2}}.
\ene
As we have already noticed, $|a_1| \leq A = \frac{1}{\sqrt{\log\frac{e}{\ep}}}$.
In particular, $|a_1| \leq 1$. Since the function 
$t \rightarrow t \log(e/t)$ is increasing in $0 \leq t \leq 1$,
$$
a_1^2 \log\frac{2}{a_1^2} \, \leq \, a_1^2 \log\frac{e}{a_1^2}
 \ \leq \ A^2 \log\frac{e}{A^2} \, = \,
\frac{1}{\log\frac{e}{\ep}}\ \bigg(1 + \log\log\frac{e}{\ep}\bigg).
$$
Taking the square root of the right-hand side, we obtain a function
which can be majorized and absorbed by the bound given in Theorem 4.1. 
As a result, we have arrived at the following consequence of this theorem.

\vskip4mm
{\bf Theorem 4.4.} {\it Assume independent random variables $X$ and $Y$
have distribution functions $F$ and $G$ with mean zero and with $\Var(X+Y) = 1$.
If\, $L(F * G,\Phi) \leq \ep < 1$, then with some 
absolute constant $C$
$$
L(F,\Phi_{v_1})\, \leq\, C\
\frac{(\log \log\frac{4}{\ep})^2}{\sqrt{\log\frac{1}{\ep}}} +
\sqrt{\,\delta_X(N)\, \log(2/\delta_X(N))}\,,
$$
where $v_1 = \sqrt{\Var(X)}$, \ $N = 1 + \sqrt{2\log(1/\ep)}$, and
$\delta_X(N) = \int_{\{|x| > N\}} x^2\,dF(x)$.
}

\vskip4mm
It seems that in general it is not enough to know that $\Var(X) \leq 1$
and $L(F * G,\Phi) \leq \ep < 1$, in order to judge the decay of
the quadratic tail function $\delta_X(T)$ as $T \rightarrow \infty$. 
So, some additional properties should be involved.
As we will see, the entropic distance perfectly suits this idea, so 
that one can start with the entropic assumption $D(X+Y) \leq \ep$. 


\section{{\bf Application of Sapogov-type results to Gaussian 
regularization}}
\setcounter{equation}{0}

\vskip2mm
In this section we consider the stability problem in Cramer's theorem
for the regularized distributions with respect to the total variation 
norm. As a basic tool, we use Theorem 2.3.

Thus, let $X$ and $Y$ be independent random variables with distribution
functions $F$ and $G$, and with variances
$\Var(X) = v_1^2$, $\Var(Y) = v_2^2$ $(v_1,v_2 > 0, \ v_1^2 + v_2^2 = 1)$,
so that $X + Y$ has variance 1. What is not important (and is assumed 
for simplicity of notations, only), let both $X$ and $Y$ have mean zero.
As we know from Theorem 2.3, the main stability result asserts that if
$
\|F * G - \Phi\| \leq \ep < 1,
$
then
$$
\|F - \Phi_{v_1}\|\, \leq\, 
\frac{Cm(v_1,\ep)}{v_1 \sqrt{\log\frac{1}{\ep}}}\,, \qquad
\|G - \Phi_{v_2}\|\, \leq\, 
\frac{Cm(v_2,\ep)}{v_2 \sqrt{\log\frac{1}{\ep}}}
$$
for some absolute constant $C$. Here, as before
$$
m(v,\ep) = 
\min\Big\{\frac{1}{\sqrt{v}},\log \log \frac{e^e}{\ep}\Big\}, \qquad
v>0, \ \ 0 < \ep \leq 1.
$$

On the other hand, such a statement -- even in the case of equal 
variances -- is no longer true for the total variation norm. 
So, it is natural to use the Gaussian regularizations
$$
X_\sigma = X + \sigma Z, \qquad Y_\sigma = Y + \sigma Z,
$$
where $Z \sim N(0,1)$ is independent of $X$ and $Y$, and where
$\sigma$ is a (small) positive parameter. 
For definiteness, we assume that $0<\sigma \leq 1$. Note that
$$
\Var(X_\sigma) = v_1^2 + \sigma^2, \quad
\Var(Y_\sigma) = v_2^2 + \sigma^2 \qquad {\rm and} \qquad
\Var(X_\sigma+Y_\sigma) = 1 + 2\sigma^2.
$$
Denote by $F_\sigma$ and $G_\sigma $ the distributions of 
$X_\sigma$ and $Y_\sigma$, respectively. Assume $X_\sigma+Y_\sigma$ 
is almost normal in the sense of the total variation norm
and hence in the Kolmogorov distance, namely,
$$
\|F_\sigma * G_\sigma - N(0,1 + 2\sigma^2)\|\, \leq\, \frac{1}{2}\
\|F_\sigma * G_\sigma - N(0,1 + 2\sigma^2)\|_{{\rm TV}}\, \leq\, \ep \leq 1.
$$
Note that $X_\sigma+Y_\sigma = (X+Y) + \sigma\sqrt{2}\, Z$
represents the Gaussian regularization of the sum $X+Y$ with parameter
$\sigma\sqrt{2}$. One may also write
$X_\sigma + Y_\sigma = X + (Y + \sigma\sqrt{2}\, Z)$, or equivalently,
$$
\frac{X_\sigma + Y_\sigma}{\sqrt{1+2\sigma^2}} = X' + Y',
\quad
\text{where}\quad
X' = \frac{X}{\sqrt{1 + 2\sigma^2}}, \quad
Y' = \frac{Y + \sigma\sqrt{2}\, Z}{\sqrt{1 + 2\sigma^2}}.
$$
Thus, we are in position to apply Theorem 2.3 to the distributions 
of the random variables $X'$ and $Y'$ with variances
$$
v_1'^2 = \frac{v_1^2}{1 + 2\sigma^2} \quad {\rm and} \quad
v_2'^2 = \frac{v_2^2 + 2\sigma^2}{1 + 2\sigma^2}. 
$$
Using $1 + 2\sigma^2 \leq 3$, it gives
$$
\|F - \Phi_{v_1}\|\, \leq\, 
\frac{Cm(v_1',\ep)}{v_1'\sqrt{\log\frac{1}{\ep}}}
\leq \frac{3C m(v_1,\ep)}{v_1 \sqrt{\log\frac{1}{\ep}}}.
$$
Now, we apply Proposition A.2.2 $b)$ to the distributions $F$ 
and $G = \Phi_{v_1}$ with $B = v_1$ and get
$$
\big\|F_\sigma - N(0,v_1^2+\sigma^2)\big\|_{{\rm TV}}
 \, \leq \, \frac{4v_1}{\sigma}\, \|F-\Phi_{v_1}\|^{1/2} 
 \, \leq \, \frac{4v_1}{\sigma}\, 
\frac{\sqrt{3C m(v_1,\ep)}}{v_1^{1/2} (\log\frac{1}{\ep})^{1/4}}.
$$
One may simplify this bound by using 
$v_1 \sqrt{m(v_1,\ep)} \leq \sqrt{v_1}$, and then
we may conclude:

\vskip4mm
{\bf Theorem 5.1.} {\it Let $F$ and $G$ be distribution functions
with mean zero and variances $v_1^2,v_2^2$, respectively, such that 
$v_1^2 + v_2^2 = 1$. Let $0 < \sigma \leq 1$. 
If the regularized distributions satisfy
$$
\frac{1}{2}\ \big\|F_\sigma * G_\sigma - N(0,1 + 2\sigma^2)\big\|_{{\rm TV}} 
\leq \ep \leq 1,
$$
then with some absolute constant $C$
\bee
\big\|F_\sigma - N(0,v_1^2+\sigma^2)\big\|_{{\rm TV}}
 \leq  
\frac{C}{\sigma}\, \bigg(\frac{1}{\log\frac{1}{\ep}}\bigg)^{1/4}, \quad 
\big\|G_\sigma - N(0,v_2^2+\sigma^2)\big\|_{{\rm TV}} 
 \leq 
\frac{C}{\sigma}\, \bigg(\frac{1}{\log\frac{1}{\ep}}\bigg)^{1/4}.
\ene
}


\section{{\bf Control of tails and entropic Chebyshev-type inequality}}
\setcounter{equation}{0}

One of our further aims is to find an entropic version of the Sapogov 
stability theorem for regularized distributions. As part of the problem,
we need to bound the quadratic tail function
$$
\delta_X(T) = \E X^2\, 1_{\{|X| \geq T\}}
$$
quantitatively in terms of the entropic distance $D(X)$.
Thus, assume a random variable $X$ has mean zero and variance $\Var(X)=1$, 
with a finite distance to the standard normal law
$$
D(X) = h(Z) - h(X) = \int_{-\infty}^\infty 
p(x)\, \log\frac{p(x)}{\varphi(x)}\,dx,
$$
where $p$ is density of $X$ and $\varphi$ is the density of $N(0,1)$. 
One can also write another representation,
$D(X) = \Ent_\gamma(f)$, where $f=\frac{p}{\varphi}$,
with respect to the standard Gaussian measure $\gamma$ on the real line.
Let us recall that the entropy functional
$$
\Ent_\mu(f) = \E_\mu f \log f - \E_\mu f \log \E_\mu f
$$
is well-defined for any measurable function $f \geq 0$ on an abstract 
probability space $(\Omega,\mu)$, where $\E_\mu$ stands for the 
expectation (integral) with respect to $\mu$.

We are going to involve a variational formula for this functional (cf. e.g. [Le]): 
For all measurable functions $f \geq 0$ and $g$ on $\Omega$, such that $\Ent_\mu(f)$ 
and $\E_\mu\, e^g$ are finite,
$$
\E_\mu f g \leq \Ent_\mu(f) + \E_\mu f \log \E_\mu\, e^g.
$$
Applying it on $\Omega = \R$ with $\mu = \gamma$ and 
$f = \frac{p}{\varphi}$, we notice that $\E_\mu f = 1$ and get that
$$
\int_{-\infty}^\infty p(x)\,g(x)\, dx\ \leq\ D(X) + 
\log \int_{-\infty}^\infty e^{g(x)}\, \varphi(x)\,dx.
$$
Take $g(x) = \frac{\alpha}{2}\, x^2\, 1_{\{|x| \geq T\}}$ with
a parameter $\alpha \in (0,1)$. Then,
\bee
\int_{-\infty}^\infty e^{g(x)}\, \varphi(x)\,dx 
 & = & \gamma[-T,T] +
\int_{\{|x| \geq T\}} e^{\frac{\alpha}{2}\, x^2}\, \varphi(x)\,dx \\
 & \hskip-40mm = &
\hskip-20mm \gamma[-T,T] +  \frac{2}{\sqrt{2\pi}}\, 
\int_T^\infty e^{-(1 - \alpha)\,x^2/2}\, dx
 \, = \,
\gamma[-T,T] +  \frac{2}{\sqrt{1-\alpha}}
\left(1 - \Phi\left(T \sqrt{1 - \alpha} \, \right)\right).
\ene
Using $\gamma[-T,T] < 1$ and the inequality $\log(1 + t) \leq t$, we obtain that
$$
\log \int_{-\infty}^\infty e^{g(x)}\, \varphi(x)\,dx\ \leq\
 \frac{2}{\sqrt{1-\alpha}}
\left(1 - \Phi\left(T \sqrt{1 - \alpha}\,\right)\right).
$$
Therefore,
$$
\frac{1}{2}\, \delta_X(T)\, \leq\, \frac{1}{\alpha}\, D(X) +
\frac{2}{\alpha \sqrt{1-\alpha}}
\left(1 - \Phi\left(T \sqrt{1 - \alpha}\,\right)\right).
$$

Now, we need to optimize the right-hand side over all $\alpha \in (0,1)$.
First, the standard bound $1 - \Phi(t) \leq \varphi(t)/t$ gives
\be
\frac{1}{2}\, \delta_X(T)\, \leq\, \frac{1}{\alpha}\, D(X) +
\frac{2}{\sqrt{2\pi}}\, 
\frac{1}{T \alpha (1 - \alpha)}\ e^{-(1 - \alpha)\, T^2/2}.
\en
Choosing just $\alpha = 1/2$, we get
$$
\frac{1}{2}\, \delta_X(T)\, \leq\, 2 D(X) +
\frac{8}{T\sqrt{2\pi}}\, e^{-T^2/4}\, \leq\, 2 D(X) + 2\, e^{-T^2/4},
$$
where the last bounds is fulfilled for $T \geq 4/\sqrt{2\pi}$.
For the remaining $T$ the obtained inequality is fulfilled 
automatically, since then $2 e^{-T^2/4} \geq 2 e^{-4/2\pi} > 1$,
while $\frac{1}{2}\,\delta_X(T) \leq \frac{1}{2}\,\E X^2 = \frac{1}{2}$.

Thus, we have proved the following:

\vskip4mm
{\bf Proposition 6.1.} {\it If $X$ is a random variable with $\E X = 0$ 
and $\Var(X) = 1$, having density $p(x)$, then for all $T>0$,
$$
\int_{\{|x| \geq T\}} x^2\,p(x)\, dx\, \leq\, 4 D(X) + 4\, e^{-T^2/4}.
$$
In particular, the above integral does not exceed $8 D(X)$ for
$T= 2\sqrt{\log^+(1/D(X))}$.
}

\vskip5mm
The choice $\alpha = 2/T^2$ in (6.1) would lead to a better asymptotic 
in $T$. Indeed, if $T \geq 2$, then $T\alpha (1-\alpha) \geq 1/T$, so
$$
\frac{1}{2}\, \delta_X(T)\, \leq\, \frac{T^2}{2}\, D(X) +
\frac{2e\,T}{\sqrt{2\pi}}\, e^{-T^2/2}\, \leq\, 
\frac{T^2}{2}\,  D(X) + 3T\, e^{-T^2/2}.
$$
Hence, we also have:

\vskip4mm
{\bf Proposition 6.2.} {\it If $X$ is a random variable with $\E X = 0$ 
and $\Var(X) = 1$, having density $p(x)$, then for all $T \geq 2$,
$$
\int_{\{|x| \geq T\}} x^2\,p(x)\, dx\, \leq\, T^2 D(X) + 6T\, e^{-T^2/2}.
$$
}

In the Gaussian case $X=Z$ this gives an asymptotically correct bound
for $T \rightarrow \infty$ (up to a factor).
Note as well that in the non-Gaussian case, from
Proposition 6.1 we obtain an entropic Chebyshev-type inequality
$$
\P\left\{|X| \geq 2\sqrt{\log(1/D(X))}\,\right\}\, \leq\, 
\frac{2 D(X)}{\log(1/D(X))} \qquad (D(X) < 1).
$$

Finally, let us give a more flexible variant of Propositions 6.1
with an arbitrary variance $B^2 = \Var(X)$ ($B > 0$), but still with 
mean zero. Applying the obtained statements to the random variable 
$X/B$ and replacing the variable $T$ with $T/B$, we then get that
$$
\frac{1}{B^2} \int_{\{|x| \geq T\}} x^2\,p(x)\, dx\, \leq\, 
4 D(X) + 4\, e^{-T^2/4B^2}.
$$

\section{{\bf Entropic control of tails for sums of independent summands}}
\setcounter{equation}{0}

\vskip2mm
We apply Proposition 6.1 in the following situation. Assume we have two independent 
random variables $X$ and $Y$ with mean zero, but perhaps with different variances 
$\Var(X)$ and $\Var(Y)$. Assume they have densities. 
The question is: Can we bound the tail functions $\delta_X$ and $\delta_Y$ 
in terms of $D(X+Y)$, rather than in terms of $D(X)$ and $D(Y)$?
In case $\Var(X+Y) = 1$, by Proposition 6.1, applied to the sum $X+Y$,
\be
\delta_{X+Y}(T) = \E\, (X+Y)^2\, 1_{\{|X+Y| \geq T\}} \leq
4\, D(X+Y) + 4\, e^{-T^2/4}.
\en
Hence, to answer the question, it would be sufficient to bound from below 
the tail functions $\delta_{X+Y}$ in terms of $\delta_X$ and $\delta_Y$.

Assume for a while that $\Var(X+Y) = 1/2$. In particular,
$\Var(Y) \leq 1/2$, and according to the usual Chebyshev's inequality, 
$\P\{Y \geq -1\} \geq  \frac{1}{2}$. Hence, for all $T \geq 0$, 
\bee
\E\, (X+Y)^2\, 1_{\{X+Y \geq T\}} & \geq &
\E\, (X+Y)^2\, 1_{\{X\geq T+1,\, Y \geq -1\}} \\
 & \geq &
\E\, (X-1)^2\, 1_{\{X\geq T+1,\, Y \geq -1\}}
 \ \geq \ 
\frac{1}{2}\, \E\, (X-1)^2\, 1_{\{X\geq T+1\}}. \\
\ene
If $X \geq T+1 \geq 4$, then clearly $(X-1)^2 \geq \frac{1}{2}\, X^2$, 
hence,
$
\E\, (X-1)^2\, 1_{\{X\geq T+1\}} \geq 
\frac{1}{2}\, \E\, X^2\, 1_{\{X\geq T+1\}}.
$
With a similar bound for the range $X \leq -(T+1)$, we get
\be
\delta_{X+Y}(T) \geq \frac{1}{4}\,\delta_{X}(T+1), \qquad T \geq 3.
\en
Now, change $T+1$ with $T$ (assuming that $T \geq 4$) and apply (7.1) 
to $\sqrt{2}\, (X+Y)$. Together with (7.2) it gives
$
\frac{1}{4}\,\delta_{\sqrt{2} X}(T) \leq 4\, 
D\big(\sqrt{2}\, (X + Y)\big) + 4\, e^{-(T-1)^2/4}.
$
But the entropic distance to the normal is invariant under rescaling 
of coordinates, i.e., $D(\sqrt{2}\, (X + Y)) = D(X + Y)$.
Since also $\delta_{\sqrt{2} X}(T) = 2\,\delta_{X}(T/\sqrt{2})$,
we obtain that
$$
\delta_X(T/\sqrt{2}) \leq 8\, D(X+Y) + 8\, e^{-(T-1)^2/4},
$$
provided that $T \geq 4$. Simplifying by 
$e^{-(T-1)^2/4} \leq e^{-T^2/8}$ (valid for $T \geq 4$), and then 
replacing $T$ with $T\sqrt{2}$, we arrive at
$$
\delta_X(T) \leq 8\, D(X+Y) + 8\, e^{-T^2/4}, \qquad T \geq 4/\sqrt{2}.
$$
Finally, to involve the values $0 \leq T \leq 4/\sqrt{2}$, just use 
$e^2 < 8$, so that the above inequality holds automatically for this 
range:\quad
$
\delta_X(T) \leq \Var(X) \leq 1 < 8\, e^{-T^2/4}.
$
Moreover, in order to allow an arbitrary variance $\Var(X+Y) = B^2$ 
($B>0$), the above estimate should be applied to $X/B\sqrt{2}$ and
$Y/B\sqrt{2}$ with $T$ replaced by $T/B\sqrt{2}$. Then it takes the form
$$
\frac{1}{2B^2}\, \delta_X(T) \leq 8\, D(X+Y) + 8\, e^{-T^2/8B^2}.
$$

We can summarize.

\vskip4mm
{\bf Proposition 7.1.} {\it Let $X$ and $Y$ be independent random 
variables with mean zero and with $\Var(X+Y) = B^2$ $(B>0)$. Assume $X$ 
has a density $p$. Then, for all $T \geq 0$,
$$
\frac{1}{B^2} \int_{\{|x| \geq T\}} x^2\,p(x)\, dx\, \leq\, 
16\,D(X+Y) + 16\, e^{-T^2/8B^2}.
$$
}


\vskip5mm
\section{{\bf Stability for L\'evy distance under entropic hypothesis}}
\setcounter{equation}{0}

\vskip2mm
Now we can return to the variant of the Chistyakov-Golinski result, 
as in Theorem 4.4. Let the independent random variables 
$X$ and $Y$ have mean zero, with $\Var(X+Y) = 1$, and denote 
by $F$ and $G$ their distribution functions. Also assume $X$ has a density $p$.
In order to control the term $\delta_X(N)$ in Theorem 4.4, we are going 
to impose the stronger condition
$$
D(X+Y) \leq 2\ep.
$$
Using Pinsker's inequality, this yields bounds 
for the total variation and Kolmogorov distances
$$
\|F*G - \Phi\|\, \leq\, \frac{1}{2}\,
\|F*G - \Phi\|_{{\rm TV}}\, \leq\,
\frac{1}{2}\, \sqrt{2 D(X+Y)}\, \leq\, \sqrt{\ep} = \ep'.
$$
Hence, the assumption of Theorem 4.4 is fulfilled, whenever $\ep < 1$.

As for the conclusion, first apply Proposition 7.1 with $B = 1$, which 
gives
$$
\delta_X(T) = \int_{\{|x| \geq T\}} x^2\,p(x)\, dx\, \leq\, 
16\,D(X+Y) + 16\, e^{-T^2/8}\, \leq\, 16\,\ep + 16\, e^{-T^2/8}.
$$
In our situation, $N = 1 + \sqrt{2\log(1/\ep')}= 1 + \sqrt{\log(1/\ep)}$, 
so, $\delta_X(N) \leq 16\,\ep + 16\, e^{-N^2/8} \leq C \ep^{1/8}$.
Thus, we arrive at:

\vskip4mm
{\bf Proposition 8.1.} {\it Let the independent random variables $X$ and $Y$ 
have mean zero, with $\Var(X+Y) = 1$, and assume that $X$ has a density with 
distribution function $F$. If $D(X+Y) \leq 2\ep < 2$, then
$$
L(F,\Phi_{v_1})\, \leq\, C\
\frac{(\log \log\frac{4}{\ep})^2}{\sqrt{\log\frac{1}{\ep}}},
$$
where $v_1 = \sqrt{\Var(X)}$ and $C$ is an absolute constant.
}

\vskip4mm
In general, in the conclusion one cannot replace the L\'evy distance 
$L(F,\Phi_{v_1})$ with the entropic destance $D(X)$. However, this is
indeed possible for regularized distributions, as we will see
in the next sections.

\section{{\bf Entropic distance and uniform deviation of densities}}
\setcounter{equation}{0}
Let $X$ and $Y$ be independent random variables with mean zero,
finite variances, and assume $X$ has a bounded density $p$. 
Our next aim is to estimate the entropic distance to the normal, 
$D(X)$, in terms of $D(X+Y)$ and the uniform deviation of $p$ above 
the normal density
$$
\Delta(X) = {\rm ess\, sup}_x\, (p(x) - \varphi_v(x)),
$$
where $v^2 = \Var(X)$ and $\varphi_v$ stands for the density of the 
normal law $N(0,v^2)$.

For a while, assume that $\Var(X) = 1$. Proposition A.3.2 gives the preliminary estimate
$$
D(X) \, \leq \, \Delta(X) \left[\sqrt{2\pi} + 2T + 
2T\, \log\left(1 + \Delta(X) \sqrt{2\pi}\, e^{T^2/2}\right)\right] +
\frac{1}{2}\, \delta_X(T),
$$
involving the quadratic tail function $\delta_X(T)$.
In the general situation one cannot say anything definite about the decay 
of this function. However, it can be bounded in terms of $D(X+Y)$
by virtue of Proposition 7.1: we know that, for all $T \geq 0$,
$$
\frac{1}{2B^2}\, \delta_X(T)\, \leq\, 8\,D(X+Y) + 8\, e^{-T^2/8B^2},
$$
where $B^2 = \Var(X+Y) = 1 + \Var(Y)$. So, combining the two estimates 
yields
\bee
D(X) & \leq & 8B^2\,D(X+Y) + 8B^2\, e^{-T^2/8B^2} \\
 & & +\ 
\Delta \left[\sqrt{2\pi} + 2T + 
2T\, \log\big(1 + \Delta \sqrt{2\pi}\, e^{T^2/2}\big)\right], \quad
{\rm where} \ \ \Delta = \Delta(X).
\ene
First assume $\Delta \leq 1$ and apply the above with
$T^2 = 8 B^2 \log \frac{1}{\Delta}$.
Then $8B^2\, e^{-T^2/8B^2} = 8B^2\,\Delta$, and putting 
$\beta = 4B^2 - 1 \geq 3$, we also have
\bee
\log\big(1 + \Delta \sqrt{2\pi}\, e^{T^2/2}\big)
 & = &
\log\left(1 + \Delta^{-\beta} \sqrt{2\pi}\,\right) = 
\beta \log\left(1 + \Delta^{-\beta} \sqrt{2\pi}\,\right)^{1/\beta}
 \\
 & < &
\beta \log\bigg(1 + \frac{(2\pi)^{1/2\beta}}{\Delta}\bigg)
 \ < \
\beta \log\bigg(1 + \frac{2}{\Delta}\bigg).
\ene
Collecting all the terms and using $B \geq 1$, we are lead to the 
estimate of the form
$$
D(X) \, \leq \, 8B^2\,D(X+Y) + CB^3\, \Delta
\log^{3/2}\bigg(2 + \frac{1}{\Delta}\bigg),
$$
where $C>0$ is an absolute constant.
It holds also in case $\Delta > 1$ in view of the
logarithmic bound of Proposition A.3.1,
$$
D(X) \leq \log\left(1 + \Delta\sqrt{2\pi}\,\right) + \frac{1}{2}.
$$
Therefore, the obtained bound holds true without any restriction on 
$\Delta$.

Now, to relax the variance assumption, assume
$\Var(X) = v_1^2$, $\Var(Y) = v_2^2$ ($v_1,v_2 > 0$), and without
loss of generality, let $\Var(X+Y) = v_1^2 + v_2^2 = 1$.
Apply the above to $X' = \frac{X}{v_1}$, $Y' = \frac{Y}{v_1}$.
Then, $B^2 = 1/v_1^2$ and $\Delta(X') = v_1\, \Delta(X)$, 
so with some absolute constant $c>0$,
$$
c\,v_1^2\, D(X) \, \leq \,
D(X+Y) + \Delta(X) \log^{3/2}\bigg(2 + \frac{1}{v_1 \Delta(X)}\bigg).
$$
As a result, we arrive at:

\vskip4mm
{\bf Proposition 9.1.} {\it Let $X,Y$ be independent random variables 
with mean zero, $\Var(X+Y)=1$, and such that
$X$ has a bounded density. Then, with some absolute constant $c>0$,
$$
c\,\Var(X)\, D(X) \, \leq \, D(X+Y) + 
\Delta(X) \log^{3/2}\bigg(2 + \frac{1}{\sqrt{\Var(X)}\, \Delta(X)}\bigg).
$$
}
Replacing the role of $X$ and $Y$, and adding the two inequalities,
we also have as corollary:

\vskip4mm
{\bf Proposition 9.2.} {\it Let $X,Y$ be independent random variables 
with mean zero and positive variances $v_1^2 = \Var(X)$, 
$v_2^2 = \Var(Y)$, such that $v_1^2 + v_2^2 = 1$, and both with 
densities. Then, with some absolute constant $c>0$,
$$
c\, (v_1^2\, D(X) + v_2^2\, D(Y)) \leq D(X+Y) + 
\Delta(X) \log^{3/2}\bigg(2 + \frac{1}{v_1 \Delta(X)}\bigg)
+\Delta(Y) \log^{3/2}\bigg(2 + \frac{1}{v_2 \Delta(Y)}\bigg).
$$
}

This inequality may be viewed as the inverse to the general property of 
the entropic distance, which we mentioned before, namely,
$v_1^2\, D(X) + v_2^2\, D(Y) \geq D(X+Y)$, 
under the normalization assumption $v_1^2 + v_2^2 = 1$.
Let us also state separately Proposition 9.1 in the particular case of
equal unit variances, keeping the explicit constant $8B^2 = 16$ 
in front of $D(X+Y)$. 

\vskip4mm
{\bf Proposition 9.3.} {\it Let $X,Y$ be independent random variables 
with mean zero and variances $\Var(X) = \Var(Y) = 1$, and such that
$X$ has a density. Then, with some absolute constant $C$
$$
D(X) \, \leq \, 16\,D(X+Y)+ C\, \Delta(X) 
\log^{3/2}\bigg(2 + \frac{1}{\Delta(X)}\bigg).
$$
}

One may simplify the right-hand side for small values of $\Delta(X)$ 
and get a slightly weaker inequality
$D(X) \, \leq \, 16\,D(X+Y) + C_\alpha\, \Delta(X)^\alpha$, 
$0 < \alpha < 1$, where the constants $C_\alpha$ depend on $\alpha$, only.
For large values of $\Delta(X)$, the above inequality holds, as well,
in view of the logarithmic bound of Proposition of A.3.1.


\section{{\bf The case of equal variances}}
\setcounter{equation}{0}

We are prepared to derive an entropic variant of Sapogov-type stability
theorem for regularized distributions. That is, we are going to estimate 
$D(X_\sigma)$ and $D(Y_\sigma)$
in terms of $D(X_\sigma + Y_\sigma)$ for two independent random variables 
$X$ and $Y$ with distribution functions $F$ and $G$, by involving a small 
``smoothing'' parameter $\sigma>0$. It will not be important whether or 
not they have densities. Since it will not be important
for the final statements, let $X$ and $Y$ have mean zero. 
Recall that, given $\sigma>0$, the regularized random variables are 
defined by $X_\sigma = X + \sigma Z$, $Y_\sigma = Y + \sigma Z$,
where $Z$ is independent of $X$ and $Y$, and has a standard normal 
density $\varphi$. The distributions of $X_\sigma,Y_\sigma$ are denoted 
$F_\sigma,G_\sigma$, with densities $p_\sigma,q_\sigma$. 

In this section, we consider the case of equal variances, say, 
$\Var(X) = \Var(Y) = 1$. Put
$$
\sigma_1 = \sqrt{1 + \sigma^2}, \qquad \sigma_2 = \sqrt{1 + 2\sigma^2}.
$$
Since $\Var(X_\sigma) = \Var(Y_\sigma) = \sigma_1^2$, the corresponding 
entropic distances are given by
$$
D(X_\sigma) = h(\sigma_1 Z) - h(X_\sigma) = 
\int_{-\infty}^\infty 
p_\sigma(x)\, \log\frac{p_\sigma(x)}{\varphi_{\sigma_1}(x)}\,dx,
$$
and similarly for $Y_\sigma$, where, as before, $\varphi_v$
represents the density of $N(0,v^2)$.
Assume that $D(X_\sigma + Y_\sigma)$ is small in the sense that
$D(X_\sigma + Y_\sigma) \leq 2\ep < 2$. According to Pinsker's inequality,
this yields bounds for the total variation and Kolmogorov distances
$$
\|F_\sigma * G_\sigma - \Phi_{\sigma_2}\|\, \leq\, \frac{1}{2}\,
\|F_\sigma * G_\sigma - \Phi_{\sigma_2}\|_{{\rm TV}}\, \leq\,
\sqrt{\ep} < 1.
$$

In the sequel, let $0 < \sigma \leq 1$. This guarantees that the ratio of 
variances of the components in the convolution
$F_\sigma * G_\sigma = F * (G * \Phi_{\sigma \sqrt{2}}\,)$ is bounded away
from zero by an absolute constant, so that we can apply Theorem 2.3.
Namely, it gives that $\|F - \Phi\| \leq C\log^{-1/2}(\frac{1}{\ep})$,
and similarly for $G$. (Note that raising $\ep$ to any positive power 
does not change the above estimate.) Applying Proposition A.2.1 $a)$, 
when one of the distributions is normal, we get
$$
\Delta(X_\sigma) \, = \,
\sup_x\, (p_\sigma(x) - \varphi_{\sigma_1}(x))\, \leq\, 
\frac{1}{\sigma}\, \|F - \Phi\| \, \leq\, 
\frac{C}{\sigma \sqrt{\log\frac{1}{\ep}}}.
$$

We are in position to apply Proposition 9.3 to the
random variables $X_\sigma/\sigma_1$, $Y_\sigma/\sigma_1$.
It gives
$$
D(X_\sigma) \, \leq \, 16\,D(X_\sigma + Y_\sigma) + C\, \Delta(X_\sigma) 
\log^{3/2}\Big(2 + \frac{1}{\Delta(X_\sigma)}\Big) \, \leq \,
32\,\ep + C'\,\frac{\log^{3/2}\big(2 + 
\sigma\sqrt{\log\frac{1}{\ep}}\,\big)}{\sigma\sqrt{\log\frac{1}{\ep}}},
$$
where $C'$ is an absolute constant. In the last expression 
the second term dominates the first one, and at this point, 
the assumption on the means may be removed. We arrive at:

\vskip4mm
{\bf Proposition 10.1.} {\it Let $X,Y$ be independent random 
variables with mean zero and variance one. Given $0 < \ep  < 1$ 
and $0 < \sigma \leq 1$, the regularized random variables 
$X_\sigma$, $Y_\sigma$ satisfy
\be
D(X_\sigma + Y_\sigma) \leq 2\ep \ \Rightarrow \
D(X_\sigma) + D(Y_\sigma) \leq C\,
\frac{\log^{3/2}\big(2 + 
\sigma\sqrt{\log\frac{1}{\ep}}\,\big)}{\sigma\sqrt{\log\frac{1}{\ep}}},
\en
where $C$ is an absolute constant.
}

\vskip4mm
This statement may be formulated equivalently by solving the above 
inequality with respect to $\ep$. The function
$u(x) = \frac{x}{\log^{3/2}(2+x)}$ is increasing in $x \geq 0$,
and, for any $a \geq 0$, 
$u(x) \leq a \Rightarrow x \leq 8\,a\log^{3/2}(2+a)$.
Hence, assuming $D(X_\sigma + Y_\sigma) \leq 1$, we obtain from (10.1) that
$$
\sigma\sqrt{\log\frac{1}{\ep}} \leq \frac{8C}{D}\,\log^{3/2}(2+C/D)
\leq \frac{C'}{D}\,\log^{3/2}(2 + 1/D)
$$
with some absolute constant $C'$, where $D = D(X_\sigma) + D(Y_\sigma)$. As a result,
$$
D(X_\sigma + Y_\sigma) \geq 
\exp\Big\{-\frac{C'^2 \log^3 (2 + 1/D)}{\sigma^2 D^2}\Big\}.
$$
Note also that this inequality is fulfilled automatically, if 
$D(X_\sigma + Y_\sigma) \geq 1$. Thus, we get:

\vskip4mm
{\bf Proposition 10.2.} {\it Let $X,Y$ be independent random variables 
with $\Var(X) = \Var(Y) = 1$. Given $0 < \sigma \leq 1$,
the regularized random variables $X_\sigma$ and $Y_\sigma$ satisfy
$$
D(X_\sigma + Y_\sigma) \geq 
\exp\Big\{-\frac{C\,\log^3 (2 + 1/D)}{\sigma^2 D^2}\Big\},
$$
where $D = D(X_\sigma) + D(Y_\sigma)$ and $C>0$ is an absolute constant.
}


\section{{\bf Proof of Theorem 1.1}}
\setcounter{equation}{0}

Now let us consider the case of arbitrary variances
$$
\Var(X) = v_1^2, \quad \Var(Y) = v_2^2 \quad (v_1,v_2 \geq 0).
$$
For normalization reasons, let $v_1^2 + v_2^2 = 1$. Then
$$
\Var(X_\sigma) = v_1^2 + \sigma^2, \quad \Var(Y_\sigma) = v_2^2 + \sigma^2, 
\quad \Var(X_\sigma+  Y_\sigma) = \sigma_2^2,
$$
where $\sigma_2 = \sqrt{1+2\sigma^2}$.
As before, we assume that both $X$ and $Y$ have mean zero, although this 
will not be important for the final conclusion.

Again, we start with the hypothesis $D(X_\sigma + Y_\sigma) \leq 2\ep < 2$
and apply Pinsker's inequality:
$$
\|F_\sigma * G_\sigma - \Phi_{\sigma_2}\|\, \leq\, \frac{1}{2}\,
\|F_\sigma * G_\sigma - \Phi_{\sigma_2}\|_{{\rm TV}}\, \leq\,
\sqrt{\ep} < 1.
$$
For $0 < \sigma \leq 1$, write $F_\sigma * G_\sigma = F * (G * \Phi_{\sigma \sqrt{2}}\,)$.
Now, the ratio of variances of the components in the convolution,
$\frac{v_1^2}{1 + 2\sigma^2}$, may not be bounded away from zero, since 
$v_1$ is allowed to be small. Hence, the application of Theorem 2.3 will only give
$
\|F - \Phi_{v_1}\| \leq \frac{Cm(v_1,\ep)}{v_1 \sqrt{\log\frac{1}{\ep}}}
$
and similarly for $G$. The appearance of $v_1$ on the right is 
however not desirable. So, it is better to involve the L\'evy distance, 
which is more appropriate in such a situation.
Consider the random variables 
$$
X' = \frac{X}{\sqrt{1+2\sigma^2}}, \quad
Y' = \frac{Y + \sigma \sqrt{2} Z}{\sqrt{1+2\sigma^2}},
$$
so that $\Var(X'+Y')=1$, and denote by $F'$, $G'$ their distribution 
functions. Since the Kolmogorov distance does not change after rescaling 
of the coordinates, we still have
$$
L(F' * G',\Phi) \, \leq \,  \|F' * G' - \Phi\|  \, = \,
\|F_\sigma * G_\sigma - \Phi_{\sigma_2}\|\, \leq\, \sqrt{\ep} < 1.
$$
In this situation, we may apply Proposition 8.1 to the couple $(F',G')$. 
It gives that
$$
L(F',\Phi_{v_1'})\, \leq\, C\ \Big(\log \log\frac{4}{\ep}\Big)^2\,
\Big(\log\frac{1}{\ep}\Big)^{-1/2}
$$
with some absolute constant $C$, where 
$
v_1' = \sqrt{\Var(X')} = \frac{v_1}{\sqrt{1+2\sigma^2}}.
$
Since $v_1' \leq v_1 \leq \sqrt{3} v_1'$, we have a similar conclusion
about the original distribution functions, i.e.
$
L(F,\Phi_{v_1}) \leq C\, (\log \log\frac{4}{\ep})^2\, (\log\frac{1}{\ep})^{-1/2}.
$
Now we use Proposition A.2.3 (applied when one of the 
distributions is normal), which for $\sigma \leq 1$ gives
$
\Delta(X_\sigma)\, \leq\, \frac{3}{2\sigma^2}\, L(F,\Phi_{v_1}),
$
and similarly for $Y$. Hence,
\be
\Delta(X_\sigma)\, \leq\, C\
\frac{(\log \log\frac{4}{\ep})^2}{\sigma^2 \sqrt{\log\frac{1}{\ep}}}, 
\qquad
\Delta(Y_\sigma)\, \leq\, C\
\frac{(\log \log\frac{4}{\ep})^2}{\sigma^2 \sqrt{\log\frac{1}{\ep}}}.
\en
We are now in a position to apply Proposition 9.2 to the random variables 
$X_\sigma' = X_\sigma/\sqrt{1+\sigma^2}$, 
$Y_\sigma' = Y_\sigma/\sqrt{1+\sigma^2}$,
which ensures that with some absolute constant $c>0$
\bee
c\, (v_1(\sigma)^2\, D(X_\sigma) + v_2(\sigma)^2\, D(Y_\sigma))
 & \leq &  D(X_\sigma + Y_\sigma) \\
 & & \hskip-27mm +  \
\Delta(X_\sigma) \log^{3/2}\Big(2 + \frac{1}{v_1(\sigma) \Delta(X_\sigma)}\Big) +
\Delta(Y_\sigma) \log^{3/2}\Big(2 + \frac{1}{v_2(\sigma) \Delta(Y_\sigma)}\Big),
\ene
where
$v_1(\sigma)^2 = \Var(X_\sigma') = \frac{v_1^2 + \sigma^2}{1 + \sigma^2}$ and
$v_2(\sigma)^2 = \Var(Y_\sigma') = \frac{v_2^2 + \sigma^2}{1 + \sigma^2}$
($v_1(\sigma), v_2(\sigma) \geq 0$).
Note that $v_1(\sigma) \geq \sigma/\sqrt{2}$. Applying the bounds in 
(11.1), we obtain that
$$
c\, (v_1(\sigma)^2\, D(X_\sigma) + v_2(\sigma)^2\, D(Y_\sigma))
 \, \leq \,  D(X_\sigma + Y_\sigma) + 
\frac{(\log \log\frac{4}{\ep})^2}{\sigma^2 \sqrt{\log\frac{1}{\ep}}} \,
\log^{3/2}\bigg(2 + 
\frac{\sigma \sqrt{\log\frac{1}{\ep}}}{(\log \log\frac{4}{\ep})^2}\bigg)
$$
with some other absolute constant $c>0$. Here, 
$D(X_\sigma + Y_\sigma) \leq 2\ep$, which is dominated by the last
expression, and we arrive at:

\vskip4mm
{\bf Proposition 11.1.} {\it Let $X,Y$ be independent random variables with total 
variance one. Given $0 < \sigma \leq 1$, if the regularized random variables 
$X_\sigma$, $Y_\sigma$ satisfy $D(X_\sigma + Y_\sigma) \leq 2\ep < 2$, then with
some absolute constant $C$
\be
\Var(X_\sigma)\, D(X_\sigma) + \Var(Y_\sigma)\, D(Y_\sigma)
 \, \leq \, 
C\,\frac{(\log \log\frac{4}{\ep})^2}{\sigma^2 \sqrt{\log\frac{1}{\ep}}} \,
\log^{3/2}\bigg(2 + 
\frac{\sigma \sqrt{\log\frac{1}{\ep}}}{(\log \log\frac{4}{\ep})^2}\bigg).
\en
}

It remains to solve this inequality with respect to $\ep$.
Denote by $D'$ the left-hand side of (11.2) and let $D = \sigma^2 D'$.
Assuming that $D(X_\sigma + Y_\sigma) < 2$ and arguing as in the
proof of Proposition 10.2, we get
$\frac{\sigma \sqrt{\log\frac{1}{\ep}}}{(\log \log\frac{4}{\ep})^2}
\leq \frac{8C}{\sigma D'}\,\log^{3/2}(2+C/D')$, hence
$\frac{\log\frac{1}{\ep}}{(\log \log\frac{4}{\ep})^4}
\leq A \equiv \frac{C'}{D^2}\,\log^3(2 + 1/D)$
with some absolute constant $C'$. The latter inequality implies
with some absolute constants
$$
\log\frac{1}{\ep} \leq C'' A\log^4(2 + A) \leq \frac{C'''}{D^2}\,
\log^7(2 + 1/D),
$$
and we arrive at the inequality of Theorem 1.1
(which holds automatically, if $D(X_\sigma + Y_\sigma) \geq 1$).


\section{{\bf Appendix I: General bounds for distances between distribution functions}}
\setcounter{equation}{0}
Here we collect a few elementary and basically known relations 
for classical metrics, introduced at the beginning of Section 2.
Let $F$ and $G$ be arbitrary distribution functions of some random 
variables $X$ and $Y$.
First of all, the L\'evy, Kolmogorov, and the total variation 
distances are connected by the chain of the inequalities
$
0 \leq L(F,G) \leq \|F - G\| \leq \frac{1}{2}\, \|F - G\|_{{\rm TV}} \leq 1.
$
As for the Kantorovich-Rubinshtein distance, there is the following 
well-known bound.

\vskip4mm
{\bf Proposition A.1.1.} {\it We have $L(F,G) \leq W_1(F,G)^{1/2}$.
}

%
%
\vskip4mm
{\bf Proposition A.1.2.} {\it If 
$\int_{-\infty}^\infty x^2\,dF(x) \leq B^2$ and 
$\int_{-\infty}^\infty x^2\,dG(x) \leq B^2$\, $(B \geq 0)$, then\\
\centerline{$a)$ \, $W_1(F,G) \leq 2L(F,G) + 4 B\, L(F,G)^{1/2}$ \; and \;
$b)$ \, $W_1(F,G) \leq 4 B\, \|F-G\|^{1/2}$.}
}

{\bf Proof.} It follows from the definition of the L\'evy distance $h = L(F,G)$ 
that, for all $x \in \R$,
$$
|F(x) - G(x)| \leq (F(x+h) - F(x)) + (G(x+h) - G(x)) + h.
$$
Integrating this inequality over a finite interval $(a,b)$, $a<b$, and using 
a general relation $\int_{-\infty}^\infty (F(x+y) - F(x))\,dx = y$ $(y \geq 0)$, we get
$$
\int_a^b |F(x) - G(x)|\,dx \, \leq \, h\, (2 + (b-a)).
$$
By Chebyshev's inequality, $\P\{X \geq x\} \leq \frac{B^2}{x^2}$ and
$\P\{X \leq -x\} \leq \frac{B^2}{x^2}$ ($x>0$), and similarly for $Y$. Hence, 
$|F(x) - G(x)| \leq \frac{B^2}{x^2}$, and for any $b>0$,
$$
\int_{\{|x|>b\}} |F(x) - G(x)|\,dx \, \leq \, 
\int_{\{|x|>b\}} \frac{B^2}{x^2}\,dx = \frac{2B^2}{b}.
$$
Using the previous estimate over the finite interval with $a=-b$,
we arrive at
$$
\int_{-\infty}^\infty |F(x) - G(x)|\,dx \, \leq \, 
2h\, (1 + b) + \frac{2B^2}{b}.
$$
This bound can be optimized over all $b>0$ by taking 
$b = B/\sqrt{h}$, and then we get the estimate in $a)$.
In case of the Kolmogorov distance, one can use similar arguments. Indeed,
$$
\int_{-b}^b |F(x) - G(x)|\,dx \, \leq \, 2hb \quad
{\rm with} \ \ h = \|F-G\|.
$$
Hence,
$
\int_{-\infty}^\infty |F(x) - G(x)|\,dx \, \leq \, 2h b + \frac{2B^2}{b}.
$
The optimal choice $b = B/\sqrt{h}$ leads to the second bound of the proposition.
\qed


\vskip5mm
\section{{\bf Appendix II: Relations for distances between regularized distributions}}
\setcounter{equation}{0}

\vskip2mm
Now, let us turn to the regularized random variables 
$X_\sigma = X + \sigma Z$, $Y_\sigma = Y + \sigma Z$,
where $\sigma>0$ is a fixed parameter and $Z \sim N(0,1)$ is a standard 
normal random variable independent of $X$ and $Y$. They have 
distribution functions
\bee
F_\sigma(x) & = & \int_{-\infty}^\infty F(x-y)\, d\Phi_\sigma(y) = 
\int_{-\infty}^{\infty} \Phi_\sigma(x-y)\, dF(y), \\
G_\sigma(x) & = & \int_{-\infty}^\infty G(x-y)\, d\Phi_\sigma(y) = 
\int_{-\infty}^{\infty} \Phi_\sigma(x-y)\, dG(y), \\
\ene
and densities
$$
p_\sigma(x) = 
\int_{-\infty}^\infty \varphi_\sigma(x-y)\, dF(y) = - \frac{1}{\sigma^2}
\int_{-\infty}^\infty F(x-y)\, y\,\varphi_\sigma(y)\,dy,
$$
$$
q_\sigma(x) = 
\int_{-\infty}^\infty \varphi_\sigma(x-y)\, dG(y) = - \frac{1}{\sigma^2}
\int_{-\infty}^\infty G(x-y)\, y\,\varphi_\sigma(y)\,dy.
$$
So, in terms of the Kolmogorov distance,
$$
|p_\sigma(x) - q_\sigma(x)|\, \leq\,
\frac{\|F - G\|}{\sigma^2} 
\int_{-\infty}^\infty |y|\,\varphi_\sigma(y)\,dy \\
 \, = \, \frac{2}{\sqrt{2\pi}}\,\frac{\|F - G\|}{\sigma}.
$$
Similarly,
$$
\int_{-\infty}^\infty |p_\sigma(x) - q_\sigma(x)|\,dx \, \leq \,
\frac{1}{\sigma^2} \int_{-\infty}^\infty |F(x) - G(x)|\,dx 
\int_{-\infty}^\infty |y|\,\varphi_\sigma(y)\,dy \, = \, \frac{2}{\sqrt{2\pi}}\,\frac{W_1(F,G)}{\sigma}.
$$
Simplifying with the help of $\frac{2}{\sqrt{2\pi}} < 1$,
let us state these bounds once more.

\vskip4mm
{\bf Proposition A.2.1.} {\it We have

\vskip2mm
$a)$ \
$\sup_x\, |p_\sigma(x) - q_\sigma(x)|\, \leq\, \frac{1}{\sigma}\,\|F - G\|$.
\qquad
$b)$ \,
$\|F_\sigma - G_\sigma\|_{{\rm TV}}\, \leq\, \frac{1}{\sigma}\, W_1(F,G)$.
}

\vskip4mm
Thus, if $F$ is close to $G$ in a weak sense, then the regularized
distributions will be closed in a much stronger sense, at least when 
$\sigma$ is not very small.
Now, applying the general Proposition A.1.2, one may replace $W_1(F,G)$ 
in part $b)$ with other metrics:

\vskip4mm
{\bf Proposition A.2.2.} {\it If 
$\int_{-\infty}^\infty x^2\,dF(x) \leq B^2$ and 
$\int_{-\infty}^\infty x^2\,dG(x) \leq B^2$\, $(B \geq 0)$, then

\vskip3mm
$a)$ \
$ \|F_\sigma - G_\sigma\|_{{\rm TV}}\, \leq\, \frac{2}{\sigma}\,
\left[L(F,G) +  2B\,L(F,G)^{1/2}\right]$;

\vskip3mm
$b)$ \ $\|F_\sigma - G_\sigma\|_{{\rm TV}}\, \leq\,\frac{4 B}{\sigma}\, \|F-G\|^{1/2}$.
}

\vskip4mm
Combining Propositions A.1.2 and A.2.1, one may bound
$\sup_x |p_\sigma(x) - q_\sigma(x)|$ in terms of the L\'evy distance
$L(F,G)$, as well. However, in order to get rid of the unnecessary 
parameter $B$, one may argue as follows. Recall that
$$
p_\sigma(x) - q_\sigma(x)\, =\, \frac{1}{\sigma^2}
\int_{-\infty}^\infty (G(x-y) - F(x-y))\, y\,\varphi_\sigma(y)\,dy.
$$
From the definition of $h = L(F,G)$, it follows that
$
|G(u) - F(u)| \leq (G(u+h) - G(u-h)) + h,
$
for all $u \in \R$, which gives
\bee
\sigma^2\,|p_\sigma(x) - q_\sigma(x)|
 & \leq & 
\int_{-\infty}^\infty 
\big(G(x-y+h) - G(x-y-h)) + h\big)\, |y|\,\varphi_\sigma(y)\,dy \\
 & \hskip-20mm \leq & \hskip-10mm
\max_y |y|\,\varphi_\sigma(y)\,
\int_{-\infty}^\infty (G(x-y+h) - G(x-y-h))\,dy + 
h \int_{-\infty}^\infty |y|\,\varphi_\sigma(y)\,dy \\
 & \hskip-20mm = & \hskip-10mm
\frac{2h}{\sqrt{2\pi e}} + h\, \frac{2\sigma}{\sqrt{2\pi}}.
\ene
Here we used the property that the function $|y|\,\varphi_\sigma(y)$ 
is maximized at $y = \pm \sigma$. Simplifying absolute factors, 
the right-hand side can be bounded by $\frac{h}{2} + \sigma h$.
We thus obtained:

\vskip4mm
{\bf Proposition A.2.3.} {\it We have
$$
\sup_x\, |p_\sigma(x) - q_\sigma(x)|\, \leq\, \frac{L(F,G)}{\sigma}\,
\bigg(1 + \frac{1}{2\sigma}\bigg).
$$
}


\section{{\bf Appendix III: Special bounds for entropic distance to the normal}}
\setcounter{equation}{0}

Let $X$ be a random variable with mean zero and variance 
$\Var(X) = v^2$ ($v > 0$) and with a bounded density $p$. 
In this section we derive bounds for the entropic distance $D(X)$
in terms of the quadratic tail function
$
\delta_X(T) = \int_{\{|x| \geq T\}} x^2\,p(x)\,dx
$
and another quantity, which is directly responsible for the closeness 
to the normal law,
$
\Delta(X) = {\rm ess\, sup}_x\, (p(x) - \varphi_v(x)).
$
As before, $\varphi_v$ stands for the density of a normally distributed
random variable $Z \sim N(0,v^2)$, and we write $\varphi$ in the 
standard case $v = 1$.
The functional $\Delta = \Delta(X)$ is homogeneous with respect 
to $X$ with power of homogeneity $-1$ in the sense that in general
$\Delta(\lambda X) = \Delta(X)/\lambda$ $(\lambda>0)$.
Hence, the functional $\Delta = \sqrt{\Var(X)}\, \Delta(X)$ 
is invariant under rescaling of the coordinates.

To relate the two quantites, $D(X)$ and $\Delta = \Delta(X)$, first write
$
p(x) \leq \varphi_v(x) + \Delta \leq \frac{1}{v\sqrt{2\pi}} + \Delta,
$
which gives $p(x) \cdot v\sqrt{2\pi} \leq 1 + \Delta\,v\sqrt{2\pi}$.
Hence,
$$
\int_{-\infty}^\infty p(x)\, \log\frac{p(x)}{\varphi_v(x)}\,dx \, = \,
\int_{-\infty}^\infty p(x)\, 
\Big(\log\big(p(x)\, v\sqrt{2\pi}\,\big) + \frac{x^2}{2v^2}\Big)\,dx  \, \leq \,
\log\Big(1 + \Delta v\sqrt{2\pi}\,\Big) + \frac{1}{2}. 
$$
Thus we have:

\vskip4mm
{\bf Proposition A.3.1.} {\it Let $X$ be a random variable with mean 
zero and variance $\Var(X) = v^2$ $(v > 0)$, having a bounded density. 
Then
$$
D(X) \leq \log\Big(1 + v\Delta(X) \,\sqrt{2\pi}\,\Big) + \frac{1}{2}.
$$
}

This estimate might be good, when both $D(X)$ and $\Delta(X)$ 
are large, but it cannot be used to see that $X$ is almost normal.
So, we need to considerably refine Proposition A.3.1 for the case, where
$\Delta(X)$ is small. For definiteness, consider the standard case $v = 1$.
Take any $T \geq 0$. Using once more the bound
$p(x)\sqrt{2\pi} \leq 1 + \Delta\sqrt{2\pi}$, where $\Delta = \Delta(X)$,
we may write
\bee
\int_{\{|x| \geq T\}} p(x)\, \log\frac{p(x)}{\varphi(x)}\,dx 
 & = &
\int_{\{|x| \geq T\}} 
p(x)\, \Big(\log\big(p(x)\sqrt{2\pi}\,\big) + \frac{x^2}{2}\Big)\, dx \\
 & \leq &
\log\big(1 + \Delta\sqrt{2\pi}\,\big) + \frac{1}{2}\, \delta_X(T)
 \, \leq \, \Delta\sqrt{2\pi} + \frac{1}{2}\, \delta_X(T).
\ene
On the last step we used $\log(1+t) \leq t$ to simplify the bound.

For $|x| \leq T$, we use
$\frac{p(x)}{\varphi(x)} \leq 1 + \frac{\Delta}{\varphi(x)}$, so that
$
\log\frac{p(x)}{\varphi(x)} \leq 
\log(1 + \frac{\Delta}{\varphi(x)}) \leq \frac{\Delta}{\varphi(x)}.
$
This gives
\bee
p(x)\, \log\frac{p(x)}{\varphi(x)}
 & \leq &
(\varphi(x) + \Delta)\,\log\Big(1 + \frac{\Delta}{\varphi(x)}\Big) \\
 & &
\hskip-20mm \leq \ \varphi(x) \log\Big(1 + \frac{\Delta}{\varphi(x)}\Big) + 
\Delta\, \log\left(1 + \Delta \sqrt{2\pi}\, e^{T^2/2}\right)
 \, \leq \,
\Delta + \Delta\, \log\big(1 + \Delta \sqrt{2\pi}\, e^{T^2/2}\big),
\ene
and after integration over $[-T,T]$
$$
\int_{\{|x| \leq T\}} p(x)\, \log\frac{p(x)}{\varphi(x)}\,dx 
 \ \leq \
2\Delta T + 2\Delta T\, \log\big(1 + \Delta \sqrt{2\pi}\, e^{T^2/2}\big).
$$
Collecting the two bounds, we arrive at:

\vskip4mm
{\bf Proposition A.3.2.} {\it Let $X$ be a random variable with mean 
zero and variance $\Var(X) = 1$, having a bounded density. 
For all $T \geq 0$,
$$
D(X) \, \leq \, \Delta(X) \left[\sqrt{2\pi} + 2T + 
2T\, \log\big(1 + \Delta(X) \sqrt{2\pi}\, e^{T^2/2}\big)\right] +
\frac{1}{2}\, \delta_X(T).
$$
}

Hence, if $\Delta(X)$ is small and $T$ is large, but not much, the 
right-hand side can be made small. When $\Delta(X) \leq \frac{1}{2}$, one may take 
$T = \sqrt{2 \log(1/\Delta(X))}$ which leads to the estimate
$$
D(X) \, \leq \, C\, \Delta(X) \sqrt{\log(1/\Delta(X))} + 
\frac{1}{2}\, \delta_X(T),
$$
where $C$ is absolute constant. If $X$ satisfies the tail condition 
$\P\{|X| \geq t\} \leq A e^{-t^2/2}$ ($t>0$), we have
$ \delta_X(T) \leq cA\, (1+T^2)\, e^{-T^2/2}$, and then
$
D(X) \, \leq \, C\, \Delta(X) \log \frac{1}{\Delta(X)},
$
where $C$ depends on the parameter $A$, only.



\begin{thebibliography}{BH3}
\itemsep=-0pt
\small


\vskip2mm
\bibitem[B-C-G1]{B-C-G1}
Bobkov, S. G., Chistyakov, G. P., G\"otze, F.
         Entropic instability of Cramer's characterization of the normal 
         law. Selected works of Willem van Zwet, 231ñ-242, Sel. Works 
         Probab. Stat., Springer, New York, 2012.

\vskip2mm
\bibitem[B-C-G2]{B-C-G2}
Bobkov, S. G., Chistyakov, G. P., G\"otze, F. Stability problems in
         Cramer-type characterization in case of i.i.d. summands.
         Theory Probab. Appl. 57 (2012), no. 4, 701ñ-723.
         
\vskip2mm
\bibitem[B-G-R-S]{B-G-R-S}
Bobkov, S. G., Gozlan, N., Roberto, C., Samson, P.-M. Bounds on the deficit 
        in the logarithmic Sobolev inequality. J. Funct. Anal. 267 (2014), 
        no. 11, 4110--4138.   

\vskip2mm
\bibitem[C-S]{C-S} 
Carlen, E. A., Soffer, A. Entropy production by block variable summation 
         and central limit theorems. Comm. Math. Phys. 
         140 (1991), no. 2, 339--371.

\vskip2mm
\bibitem[C]{C}
Chistyakov, G. P. The sharpness of the estimates in theorems on the 
         stability of decompositions of a normal distribution and 
         a Poisson distribution. (Russian) Teor. Funkcii Funkcional. 
         Anal. i Prilozhen. Vyp. 26 (1976), 119ñ-128, iii. 

\vskip2mm
\bibitem[C-G]{C-G} 
Chistyakov, G. P., Golinskii, L. B. Order-sharp estimates for the 
         stability of decompositions of the normal distribution in the 
         Levy metric. (Russian) Translated in J. Math. Sci. 72 (1994), 
         no. 1, 2848--2871. Stability problems for stochastic models
         (Russian) (Moscow, 1991), 16--40, Vsesoyuz. Nauchno-Issled. 
         Inst. Sistem. Issled., Moscow, 1991. 

\vskip2mm
\bibitem[Cr]{Cr}       
Cram\'er, H. Ueber eine Eigenschaft der Normalen Verteilungsfunktion.
         Math. Zeitschrift, Bd. 41 (1936), 405--414.

         
\vskip2mm  
\bibitem[D-C-T]{D-C-T} 
Dembo, A., Cover, T. M., Thomas, J. A. Information-theoretic inequalities. 
         IEEE Trans. Inform. Theory, 37 (1991), no. 6, 1501--1518.

\vskip2mm  
\bibitem[F-M-P1]{F-M-P1} 
Figalli, A., Maggi, F., Pratelli, A. A refined Brunn--Minkowski inequality 
         for convex sets. Ann. Inst. H. Poincar\'e Anal. 
         Non Lin\'eaire 26 (2009), no. 6, 2511--2519. 

\vskip2mm  
\bibitem[F-M-P2]{F-M-P2} 
Figalli, A., Maggi, F., Pratelli, A. Sharp stability theorems for the 
        anisotropic Sobolev and log-Sobolev inequalities on functions of 
        bounded variation. Adv. Math. 242 (2013), 80--101. 

\vskip2mm
\bibitem[Le]{Le}
Ledoux, M. Concentration of measure and logarithmic Sobolev inequalities. 
         Seminaire de Probabilites, XXXIII, 120--216, Lecture Notes in Math., 
         1709, Springer, Berlin, 1999.

\vskip2mm
\bibitem[L1]{L1} 
Linnik, Yu. V. A remark on Cramer's theorem on the decomposition 
         of the normal law. (Russian) Teor. Veroyatnost. i Primenen. 
         1 (1956), 479--480.

\vskip2mm
\bibitem[L2]{L2} 
Linnik, Yu. V. General theorems on the factorization of infinitely
         divisible laws. III. Sufficient conditions (countable bounded
         Poisson spectrum; unbounded spectrum; ``stability''). 
         Theor. Probability Appl. 4 (1959), 142--163.

\vskip2mm
\bibitem[L-O]{L-O} 
Linnik, Yu. V., Ostrovskii, I. V. Decompositions of random variables and 
         random vectors. (Russian) Izd. Nauka, Moscow, 1972, 479 pp. 
         Translated from the Russian: Translations of Mathematical
         Monographs, Vol. 48. American Math. Society, Providence, 
         R. I., 1977, ix+380 pp. 

\vskip2mm 
\bibitem[M]{M}  
Maloshevskii, S. G. Unimprovability of N. A. Sapogov's result in 
         the stability problem of H. Cramer's theorem. (Russian) 
         Teor. Verojatnost. i Primenen, 13 (1968) 522--525. 

\vskip2mm
\bibitem[MK]{MK} 
McKean, H. P., Jr. Speed of approach to equilibrium for Kac's caricature of 
         a Maxwellian gas. Arch. Rational Mech. Anal. 21 (1966), 343--367.


\vskip2mm 
\bibitem[S1]{S1}
Sapogov, N. A. The stability problem for a theorem of Cram\'er. (Russian)  
         Izvestiya Akad. Nauk SSSR. Ser. Mat. 15 (1951), 205--218. 

\vskip2mm 
\bibitem[S2]{S2}
Sapogov, N. A. The stability problem for a theorem of Cram\'er. (Russian)  
         Izvestiya Akad. Nauk SSSR. Ser. Mat. 15 (1951), 205--218. 

\vskip2mm 
\bibitem[S3]{S3}
Sapogov, N. A. The problem of stability for a theorem of Cram\'er. 
         (Russian) Vestnik Leningrad. Univ. 10 (1955), no. 11, 61--64.         

\vskip2mm 
\bibitem[Seg]{Seg}
Segal, A. Remark on stability of Brunn--Minkowski and isoperimetric inequalities
for convex bodies, in Geometric Aspects of Functional Analysis, in: Lecture Notes in Mathematics,
2050 (2012), 381--391.            
         
\vskip2mm 
\bibitem[Se]{Se}
Senatov, V. V. Refinement of estimates of stability for a theorem of 
         H. Cram\'er. (Russian) Continuity and stability in problems of 
         probability theory and mathematical statistics. Zap. Naucn. Sem. 
         Leningrad. Otdel. Mat. Inst. Steklov. (LOMI) 61 (1976), 125--134.

\vskip2mm 
\bibitem[Sh1]{Sh1}
Shiganov, I. S. Some estimates connected with the stability of H. Cramer's 
         theorem. (Russian) Studies in mathematical statistics, 3. 
         Zap. Nauchn. Sem. Leningrad. Otdel. Mat. Inst. Steklov. 
         (LOMI) 87 (1979), 187ñ-195.

\vskip2mm 
\bibitem[Sh2]{Sh2}
Shiganov, I. S. On stability estimates of Cramer's theorem. Stability
         problems for stochastic models (Varna, 1985), 178ñ-181, 
         Lecture Notes in Math., 1233, Springer, Berlin, 1987. 


\vskip2mm 
\bibitem[Z1]{Z1}
Zolotarev, V. M. On the problem of the stability of the decomposition 
         of the normal law into components. (Russian) 
         Teor. Verojatnost. i Primenen., 13 (1968), 738--742.

\vskip2mm 
\bibitem[Z2]{Z2}
Zolotarev, V. M. Estimates of the difference between distributions in the 
Levy metric. (Russian) Collection of articles dedicated to Academician Ivan 
Matveevich Vinogradov on his eightieth birthday, I. Trudy Mat. Inst. Steklov. 
112 (1971), 224--231, 388.

\end{thebibliography}
\end{document}